\documentclass[11pt,A4]{article}
\usepackage[top=1in, bottom=1in, left=1.2in, right=1.2in]{geometry}
\usepackage{graphicx}
\usepackage{here}
\usepackage{amssymb}
\usepackage{amsmath}
\usepackage{amsthm}
\usepackage{epstopdf}
\usepackage{comment}
\newtheorem{Theorem}{Theorem}[section]
\newtheorem{Lemma}[Theorem]{Lemma}
\newtheorem{Corollary}[Theorem]{Corollary}
\newtheorem{Property}[Theorem]{Property}
\newcommand{\Proof}{\noindent{\bf Proof}\quad}
\theoremstyle{definition}

\newtheorem{Definition}[Theorem]{Definition}
\newtheorem{Example}[Theorem]{Example}
\newtheorem{Construction}[Theorem]{Construction}

\title{\bf Combinatorial Designs for Deep Learning}
\author{${}^a$Shoko Chisaki,  ${}^b$Ryoh Fuji-Hara and ${}^a$Nobuko Miyamoto\\
\\
${}^a$Department of Information Science, Tokyo University of Science\\
${}^b$Faculty of Engineering, Information and Systems, University of Tsukuba}
\date{}    

\begin{document}

\maketitle

\begin{abstract}
Deep learning is a  machine learning methodology using  multi-layer neural network.
A multi-layer neural network can be regarded  as a chain of complete bipartite  graphs.
The nodes of the first partita is the input layer and the last  is the output layer.
The edges of a bipartite graph function as  weights which are represented as a matrix.
The values of  $i$-th partita  are computed by multiplication of the weight matrix and values of $(i-1)$-th partita.
Using mass training and teacher data, the weight parameters are estimated little by little.
 Overfitting (or Overlearning) refers to a model that models the “training data” too well.
 It then becomes difficult for the model to generalize to new data which were not in the training set.
 The most popular method to  avoid overfitting is called {\it dropout}.
Dropout deletes  a random sample of activations (nodes) to zero during the training process.
A random sample  of nodes causes  more irregular frequency of dropout edges.
There is a similar sampling concept in the area of  design of experiments.
We propose a combinatorial design on dropout nodes  from each partita which balances the frequency of edges.
 We  analyze  and construct  such  designs in this paper. 
\\
\\
{\bf Keywords. } Deep learning, Dropout, Split-block design, Dropout design\\
\\
{\bf AMS classification. } 05B05,   68T05,  94C30
\end{abstract}

\section{Deep learning and Overfitting Problem}
The structure of the neural network is used for many methods of deep learning, and the model of deep learning from this background is also called a {\it deep neural network}.
Usually the expression ``deep'' refers to the number of hidden layers in the neural network. In the conventional neural network, the number of hidden layers was at most 2 or 3, but the deep neural network could have 150 hidden layers.
The deep learning model learns using large labeled data and the structure of the neural network.
This model allows us to learn feature quantities directly from the data and eliminates the need for manual feature extraction.
%By doing this, we can learn feature quantities directly from the data, and manual feature extraction is no longer necessary as it used to be.

\begin{figure}[H]
\centering
\includegraphics[width=4.5in,clip]{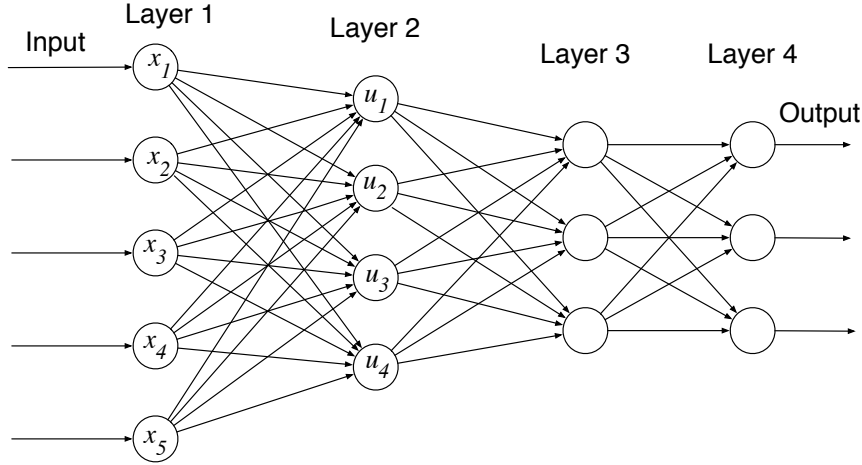}
\caption{A multi-layer neural network }
\label{MLNN}
\end{figure}

Neural networks  consist of  a series  of interconnected nodes called layers.  Networks can have tens or hundreds of hidden layers.
Consider a multi-layered neural network as shown in Figure \ref{MLNN}.
 Layer 1 is called the {\it input layer}, layers 2 and 3 are {\it internal layers or hidden layers}, and layer 4 is called the {\it output layer}.

Each connection between neurons (nodes)  is associated with a weight $w_{ij}$.
This weight dictates the importance of the input value.
The initial weights are set randomly.
Input values and hidden values are denoted as a vertical vectors $\textbf{x} = [ x_1, x_2,..., x_l ]$ and $\textbf{u}=[ u_1,u_2,...,u_m ]$, respectively, and   $\textbf{W}$ is a $m \times l $ matrix for weights.
The values of the hidden layer are calculated from the input layer as follows:
$$ \textbf{u} = \sigma \Big(  \textbf{W x} + \textbf{b} \Big),$$
where $\textbf{b}$ is  the shared value for the bias  and $\sigma$ is the neural activation function (a sigmoid function).
At each stage of the layers, the values of the next layer are calculated in the same way.
Let $\mathbf{W}^{(t)}$  and  $\mathbf{b}^{(t)}$ be the weight matrix and the bias vector of $t$-th stage, respectively.  We denote  the final  result  as $ \mathbf{y( x \,  ; W^{(1)},...,W^{(L)}, b^{(1)},...,b^{(L)} )}$ or simply  $ \mathbf{y( x \,  ; w)}$, where $\mathbf{w}$  is the vector of all weights.
Let $\mathbf{d} \ (\mathbf{d}_i)$ denote the teacher data corresponding to the input data $\mathbf{x}\  (\mathbf{x}_i)$.
%We evaluate by squared error of the difference between the teacher data and the result of the neural network.
%$$\|  \mathbf{d} - \mathbf{y(x \, ; w)} \| ^2$$

Let  $\mathbf{(x_1,d_1), (x_2,d_2),\dots, (x_N,d_N)}$  be a set of pairs of input data and teacher data． Consider  the following formula:
$$E( \mathbf{w} ) = \frac{1}{2}  \sum_{i=1}^N \|  \mathbf{d}_i - \mathbf{y(x_i \, ; w)} \|^2 .$$
We would like to choose the weights of $\mathbf{w}$ for each set of $N$ pairs of data so that $E( \mathbf{w} ) $ is minimized.
%Then, the weights of $\mathbf{w}$ are  chosen for each set of $N$ pairs of data so that  $E( \mathbf{w} ) $ is minimize.
$E(\mathbf{w})$ is called  {\it training error}.
%
%\subsection{Overfitting}

\begin{figure}[H]
\centerline {
\includegraphics[width=5in,clip]{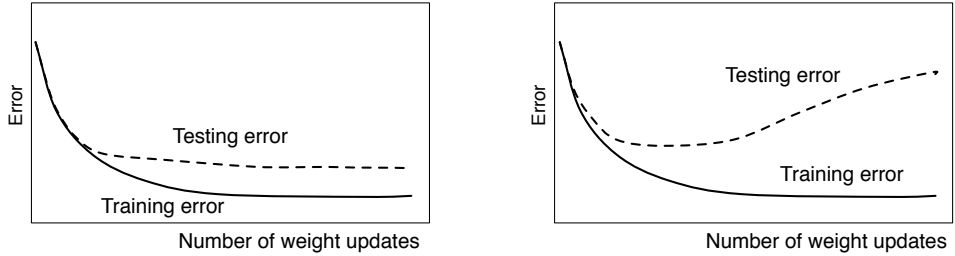}
} \caption[]{Overfitting problem }
\label{Overfitting}
\end{figure}

The real purpose of learning is to make a correct estimate for the ``unknown'' sample, which should be given from the current data, not on the given training data.
A model that fits very well (too well) for training data but not good for general data is called {\it overfitting} or {\it overlearning}.
Overfitting happens when a model learns the detail and noise in the training data to the extent that it negatively impacts the performance of the model on new data.
Therefore, a sample set different from the training data is prepared as general data, and  the error calculated by the same method as the training error is called  {\it test error}.
Training error monotonously decreases as training progresses. Ideally, the test error also decreases accordingly.
As shown in the Figure \ref{Overfitting} on the right, when the test error increases with weight update, it can be said that over learning is occurring, T. Okatani (2015) \cite{Okatani2015}.

% Dropout method

As a method to prevent overlearning, a kind of sparsity approach  called  {\it dropout} was proposed by Srivastava et al. (2014)\cite{JMLR2014} .
 In this method, neurons (nodes)  of a multilayered neural network are randomly selected and learned.
At each training stage, individual nodes are either dropped out of the net with probability $1-p$ or kept with probability $p$, so that a reduced network is left; incoming and outgoing edges to a dropped out node are also removed.  This method is widely used at present because it has good experiment results in many cases.
However, a random sample  of nodes in two layers causes  more irregular frequency of dropout edges (weights).
Let $X$ and $U$   be  random variables for  how many times node $x$ and $u$ in layer 1 and layer 2 are  selected within $n$ trials, respectively.
 Let   $V(X)$ and $V(U)$ be   variances of $X$ and $U$.  Suppose two  random variables $X$ and $U$ are converted to $Z=a X + b U$.  Then the variance  of the number of times  edge $(x,u)$ is selected can be expressed in the  form of $V(Z)=a^2V(X)+ b^2V(U)$, which implies that the edges (weights) are chosen to be more imbalance.

\begin{figure}[H]
\centerline {
\includegraphics[width=4.5in,clip]{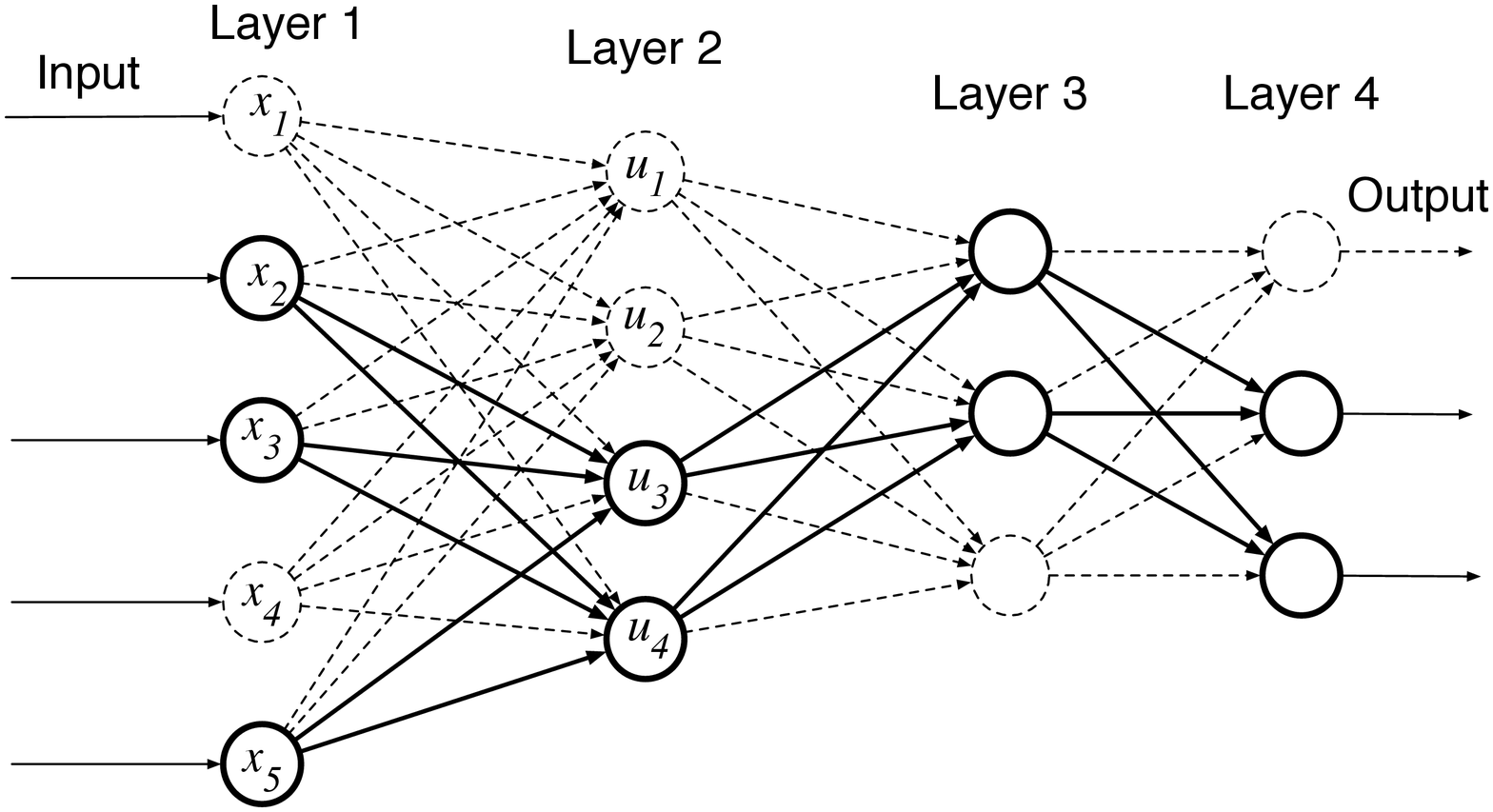}
} \caption[]{Dropout }
\end{figure}
%%%%%%%%%%%%%%%%%%%%%%%%%
\section{A statistical viewpoint}

In the first half of the 20th century,  R. A. Fisher thought that it was impossible to exclude all the factors that influence the experimental results,
and therefore the experimental results  were obliged to include fluctuations due to those influences.
On the premise of it, he thought about how to conduct experiments and lead conclusions among them.
R. A. Fisher (1935)  \cite{Fisher1935} founded an academic and practical  field called {\it the design of experiments}.
In the field of experimental design, he proposed that it is better to acquire data in a balanced manner rather than collecting data randomly for factors or treatments.
%%%%%%%%%%%%%%%%%%%%%%%

%We have discussed {\it dropout designs}  %%%%
%which satisfy  the balance conditions  instead of selecting  dropout subsets at random.
 %we show that {\it dropout designs} %%%%
%are  better than random method for weight estimation.
Here, from a statistical point of view,
we discuss on random method for weight estimation.
For example, consider a small model with two layers of 7 nodes in the input layer and 4 nodes in the hidden layer.
%, and the following {\it dropout design of type (1,1)} to use:
%$$\begin{array}{lll}
%\{x_1, x_2\ | \ {\bf u_1,u_2}\}, &\{x_3, x_5, x_6\ | \ {\bf u_1,u_2} \},& \{x_4\ | \  {\bf u_1,u_2}\},\\
%\{x_1, x_2\ | \ {\bf u_3}\}, &\{x_3, x_5, x_6 \  | \ {\bf u_3} \},& \{x_4\ | \  {\bf u_3}\}.\\
%\end{array}$$
%
%
%%%%%%%%%%%%%%%%%%%% Figures 3
 When we get input data $(x_1,x_2,...,x_7)$, the values $({\bf u_1, u_2,u_3,u_4})$ in the hidden layer are basically  determined by the following  computation:
%%%%%%%%%%% Equation
$$
\left(  \begin{array}{c} \textbf{u}_1\\ \textbf{u}_2\\ \textbf{u}_3\\ \textbf{u}_4\end{array}\right)
=
\left(  \begin{array}{ccccccc}
w_{11}& w_{12} & w_{13}& w_{14}&w_{15}&w_{16}&w_{17}\\
w_{21} &  w_{22} & w_{23}& w_{24}&w_{25}&w_{26}&w_{27}\\
w_{31} &  w_{32} & w_{33}& w_{34}&w_{35}&w_{36}&w_{37}\\
w_{41} &  w_{42} & w_{43}& w_{44}&w_{45}&w_{46}&w_{47}\\
\end{array}\right)
\left(  \begin{array}{c}  x_1 \\ x_2\\x_3\\x_4\\x_5\\x_6 \\x_7\end{array}\right).
$$
%The entries of the matrix  $\big( w_{ij} \big)$ are weights (parameters, coefficients) to estimate in the sense of statistics.
%But the  parameters are  too many comparing with data variables (explanatory variables) $x_1,x_2,...,x_6$.
In the sense of statistics, the weights $w_{ij}$ are coefficients to estimate and data variables $x_i$ are explanatory variables.
Now we focus on ${\bf u_1}$ only.  Then the equation becomes a typical linear regression model as follows:
$$
{\bf u_1} = w_{11} x_1 + w_{12} x_2 +  w_{13}x_3 +  w_{14} x_4 + w_{15}x_5 + w_{16}x_6 + w_{17}x_7 + \epsilon_1,
$$
where $\epsilon_1$ is a random variable of error.
Let $x_i^{(j)}$ be the $j$-th input data of  the variable $x_i$, $\epsilon_1^{(j)}$ be the  error for the $j$-th data and ${\bf u_1}^{(j)}$  be the value in the hidden layer determined by $j$-th input data.
%%%%%%%%%%%%%%%%%%%%%%%%%%%%%%%%%%%% Re-edit from here, Aug. 2019
%${\bf u_1}$ is obtained from the three blocks,
%$\{x_1, x_2\ | \ {\bf u_1,u_2}\}$, $\{x_3, x_5, x_6\ | \ {\bf u_1,u_2} \}$, $\{x_4\ | \  {\bf u_1,u_2}\}$,
%that is, $\mathcal{B}_2( \{ {\bf u_1}\}; V_1) = \{ \{ x_1, x_2\}, \{x_3, x_5, x_6\}, \{x_4\}\}$, where $V_1=\{  x_1,  x_2,\ldots, x_6 \}$.
%%%
Suppose that weights are randomly selected and dropped out from $w_{11}^{(j)}, w_{12}^{(j)}, ..., w_{17}^{(j)}$ every input $j=1,2,...$ .

%blocks are changed  for each $j$-th input data, $j=1,2,3,...,$  and  the block set is used repeatedly.
For example, $\textbf{u}_1^{(1)}, \textbf{u}_1^{(2)}, \ldots $ are computed by the following equation:
%\begin{center}
$$
\left( \begin{array}{c}
\textbf{u}_1^{(1)}\\
\textbf{u}_1^{(2)}\\
\textbf{u}_1^{(3)} \\
\vdots
\end{array} \right)
=
\left[ \begin{array}{lllllll}
x_1^{(1)} & x_2^{(1)} & & x_4^{(1)}&& &\\
& & x_3^{(2)} & & x_5^{(2)}& x_6^{(2)}&\\
& x_2^{(3)} &&  x_4^{(3)} & & &x_7^{(3)}\\
&&& \vdots &&&
\end{array} \right]
\left(  \begin{array}{c}  w_{11} \\ w_{12}\\w_{13}\\w_{14}\\w_{15}\\w_{16}\\w_{17} \end{array}\right) +
\left(  \begin{array}{c}
\epsilon_1^{(1)} \\
\epsilon_1^{(2)}\\
\epsilon_1^{(3)}\\
 \vdots
 \end{array}\right)
$$
%\end{center}
This is a regression model with sparse data.   The information about  variables to be saved or dropped can be expressed in the following (0,1)-matrix, called an {\it incidence matrix} or a {\it design matrix}.

%Consider the incidence matrix $X$ of the blocks $\mathcal{B}_2( \{ {\bf u_1}\}; V_1)$ by the input data $V_1$, which is

$$
X =
\left[ \begin{array}{lllllll}
1 & 1& & 1& &&\\
& &1 & & 1& 1&\\
&1 & & 1 & & &1\\
&&& \vdots &&&
\end{array} \right].
$$
\\

%\begin{figure}[H]
%\centerline {
%\includegraphics[width=2.5in,clip]{Statistics}
%} \caption[]{$\mathcal{B}_2( \{ {\bf u_1}\}; V_1)$}
%\end{figure}
When we estimate the weights  $\hat{w}_{11},    \hat{w}_{12}, \dots,  \hat{w}_{17}$ by the
regression method, it is known that they are unbiased estimations, that is, $E(\hat{w}_{1j} )= w_{1j}$.
However the precisions of estimations $V(\hat{w}_{1j} )$ are very much depend on  the patterns  of the incidence matrices $X$.
$X^t X$ is called the {\it information matrix} in Statistics.
In order to compare  the precision of estimations, we may use the determinants of information matrices.
The larger determinant of the information matrix $X^t X$ is, the better precision of the estimation (smaller variance of the estimation) is, see Nagao and Kuriki (2006) \cite{Kuriki2006}.
We had experiments for four types of random incidence matrices $X$ of size $21 \times 7$ with 63 random ones.
The four types of $X$ are in the following:
\begin{description}
\item[(1) ]completely random 63 ones in $X$,
\item[(2) ]random 63 ones  in $X$ such that  each column has exactly 9 ones,
\item[(3) ]random 63 ones  in $X$ such that each row has exactly 3 ones,
\item[(4) ]random 63 ones  in $X$ such that each column has exactly  9 ones and each row has exactly 3 ones.
\end{description}
For each type of the above, we generated  500 random incidence matrices and computed the determinants for each  information matrix, $X^tX$'s. The following chart shows  a comparison of the distributions of 500 determinants of four types by  a Box plot.
%how good the  estimation from sparse data  depends on the pattern  of the incidence matrix $X$.
\begin{figure}[H]
\centerline {
\includegraphics[width=4in,clip]{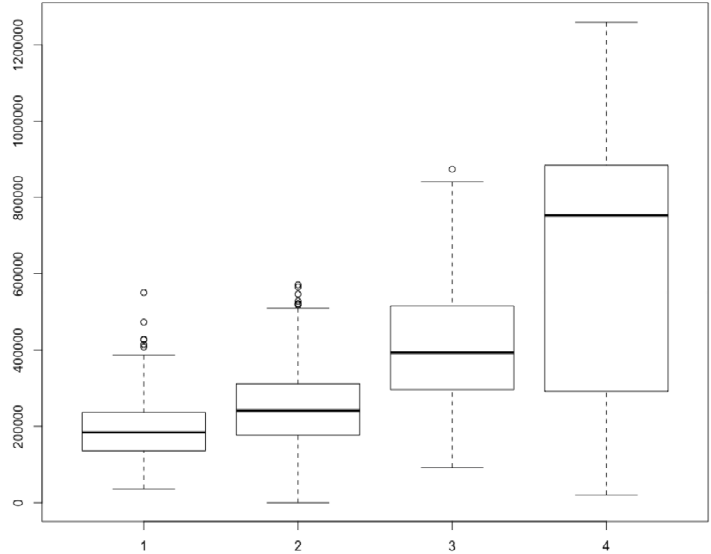}
} \caption[]{Determinant distributions of four types}
\end{figure}
In the experience of the  type (4), the maximum determinant was 1,259,712, and the information matrix  at the time was in the form  $X^t X = 6 I + 3 J$.

There are some criteria for precision of estimation, called E-optimality, A-optimality, D-optimality, (M, S)-optimality.
%E-optimality means the minimum of nonzero eigenvalues of information matrices $X^T X$,  where $X^T$ is the transpose of $X$.
If  all optimal criteria are satisfied then it is called  {\it  universally optimal}.
The next  is well known result in Statistics, see J. Kiefer (1975) \cite{Kiefer1975} .

\begin{Theorem}
\label{Kuriki}
Suppose   the number of ones in  a design matrix $X$ is a constant.
If its information matrix can be represented as
$$ X^T X = \alpha I + \beta J,$$
then the design $X$ is universally optimal,
where $\alpha, \beta$ are integers and $J$ is the all one matrix.
\end{Theorem}

Let $V=\{1,2,3,..., 7\}$, then each row of $X$ can be seen as a subset  $B_i,  i=1,..., 21$, of $V$.
Let $\mathbf{B}=\{B_1,B_2,..,B_{21}\}$. This system $(V,\mathbf{B})$ satisfies the following combinatorial conditions:
\begin{description}
\item[i. ] every element of $V$ appears in exactly $ r\, (=9)$ subsets of $\mathbf{B}$
\item[ii.] every distinct pair of $V$ appears simultaneously in exactly $\lambda \, (=3)$ subsets of $\mathbf{B}$
\end{description}
Such system is called a $(r,\lambda)${\it -design} or a {\it regular pairwise balanced design}.

%A regular $t$-wise balanced design  ($t\ge 2$) is also a regular $2$-wise balanced design. Suppose that $r$ is the number of times each point appears in the block set and $\lambda$ is the concurrence number of the regular $2$-wise balanced design.
%If  $X$ is an incidence matrix of  a  regular $t$-wise balanced design, then its information matrix is
%$$ X^T X = (r-\lambda) I + \lambda J.$$
%In  a dropout design $(V_1,V_2,...,V_n\, ; \mathcal{B} )$ of type $(d_1,d_2,...,d_t)$, consider $(V_i, V_{i+1}\ ; \mathcal{B}|_{V_i V_{i+1}} )$.
%$\mathcal{B}_{i+1}( \{ u_1\}; V_i)$, $u_1 \in V_{i+1}$,   is a balanced  $t$-wise balanced design. If  $t \ge 2$,  the incidence matrix of $\mathcal{B}_{i+1}( \{ u_1\}; V_i)$ satisfies Theorem \ref{Kuriki}.

%%%%%%%%%%%%%%%%%%%%%%%%%%%%%

\section{ Split-block designs and  related combinatorial designs }

In agricultural field experiments, sometimes similar methods to dropout are considered.
Let us consider a two-factor experiment in which a factor {\bf A} occurs at $s$ levels, $A_1, A_2, \dots, A_s$ (called {\it treatments}) and  the second factor {\bf B} occurs at $t$ levels $B_1,B_2, \dots, B_t$.
The experimenters have to obtain experiment data $y_{i,j}$ for all treatment combination
$(A_i, B_j)$ for $i=1,2,\dots,s, \ j=1,2,\dots,t$. Usually, for each treatment combination,  the experiments are done  repeatedly $\lambda$ times to estimate the treatment effects more accurately.
%Totally  $\lambda st$ data are obtained.
This gives a total of $\lambda st$ data points. Now suppose the experiments can be done  simultaneously. It is natural to break down
%Suppose we have now some experimenters and they can do the experiments simultaneously. Then we may naturally have  a thought to break down
 the two-factor experiment design to some smaller two-factor experiments called {\it blocks}.
For example,  $s=t=4$ , a $4 \times 4$ two-factorial experiment is broken down to four $2 \times 2$ two-factorial experiments (blocks) as follows:
$$\begin{array}{ll}
E_1 = \{ A_1, A_2 \ |  \ B_1, B_2\} & E_2=\{ A_3, A_4 \ | \  B_3, B_4\} \\
E_3 = \{ A_1, A_2 \ | \  B_3, B_4\} &  E_4=\{ A_3, A_4\ | \  B_1, B_2\}
\end{array}$$
Here, every treatment combination $(A_i,B_j), \  i=1,...,4,\ j=1,...,4$,  occurs exactly once in the four $2 \times 2$ two-factorial experiments, that is, $\lambda=1$.    The set of experiments above is called a {\it split-block design}, and it may have  {\it `incomplete'} or  {\it `balanced incomplete'}  as the prefix.
Now, we define a split-block design mathematically.
\begin{Definition}[Split-block design]
\label{split-block}
Let $V_1$ and $V_2$ be mutually disjoint point sets and the block set  be a collection of subsets of the points consisting of $k_1$ points from $V_1$ and $k_2$ points from $V_2$:
$$\mathcal{B} = \{ \,\{C_1 |  C_2\} \, \mid \, C_1\subset V_1, C_2 \subset V_2 , \,  |C_1|=k_1, |C_2|=k_2\}.$$
If,  for any $d_1$ and $ d_2$ points from $V_1$ and $V_2$, respectively, there exist exactly $\lambda$ blocks containing the points, then
the design $\mathcal{D}=( V_1,V_2  ;  \mathcal{B})$    is called a {\it split-block design of type $(d_1,d_2)$}.  $\lambda$ is said to be the {\it concurrence number} of the design.
\end{Definition}
$C_1$ and $C_2$ are said to be  the 1st and 2nd {\it sub-block}, respectively.  $\{C_1 | C_2\}$ is sometimes called a {\it block} or a {\it super-block}.

Designs similar to  split-block   called {\it split-plot designs} or {\it block designs with nested rows and columns}  are discussed  in 1980's.   I. Mejza (1987) \cite{MejzaI1987} first defined  the split-block designs as a development model of  split-plot design.  F. Hering and S. Mejza (1997) \cite{Hering1997} show analysis and constructions of split-block designs in more detail.

Let $b$ be the number of super-blocks of the design  $\mathcal{D}$. Let $|V_i|=v_i$ and $|C_i|=k_i$.  Then it is easy to see the following properties:
\begin{Property}[K. Ozawa et al., 2000  \cite{ozawa2000}]
$$ b k_1 k_2 = \lambda v_1 v_2 $$
\end{Property}
\begin{Property}
\label{regular}
If $\mathcal{D}$ is a split-block design of type $(d_1,d_2)$, then $\mathcal{D}$ is also a split-block design of type $(g_1,g_2)$  for any  $0\le g_1\le d_1$ and $0 \le g_2\le d_2$,  $g_1+g_2 \ge 1$.
\end{Property}
Now we show an easy construction of split-block designs.
Let $V$ be a finite set of $v$ points and $\mathcal{B}$ a collection of $k$-subsets (called blocks)  of $V$.
The pair $ (V,  \mathcal{B})$  is called a  $t$-($v,k,\lambda$) {\it design}  if every $t$-subset of $V$ appears exactly $\lambda$ times in the block set $\mathcal{B}$ , $t\ge 1$.
$\mathcal{B}_1\times \mathcal{B}_2$ is the direct product of the block sets $\mathcal{B}_1$ and $\mathcal{B}_2$:
%whose block consists of every pair of blocks from  $\mathcal{B}_1$ and $\mathcal{B}_2$:
$$ \mathcal{B}_1\times \mathcal{B}_2 = \{ \, \{C_1 | C_2\} \, \mid \, \mbox{ for all } C_1\in \mathcal{B}_1, \, C_2 \in \mathcal{B}_2 \}.$$
\begin{Construction}
\label{trivialD}
Let $(V_1,\mathcal{B}_1)$ and $(V_1,\mathcal{B}_1)$ be  $t_1$-$(v_1,k_1,\lambda_1)$ design  and $t_2$-$(v_2,k_2,\lambda_2)$ design, respectively.
$(V_1,V_2; \mathcal{B}_1\times \mathcal{B}_2)$ is a split-block design of type $(t_1,t_2)$.
The concurrence number is $\lambda_1\lambda_2$.
\end{Construction}
Let us call this design  a {\it trivial split-block design}.  The trivial split-block designs have a  bad property.
Let $b_1$ and $b_2$ be the number of blocks of the  block sets $\mathcal{B}_1$ and $ \mathcal{B}_2$, respectively.
Every  sub-block of $\mathcal{B}_1\times \mathcal{B}_2$ has $b_2$ or $b_1$ copies  in the block set.
This causes a  decrease in the variation of blocks.

\begin{Example}
Let $V_1=\{0,1,\dots ,8\} , V_2=\{ {\bf 0,1,\dots, 5} \}$.  The following is a trivial split-block design of type $(1,1)$.
The block set is  $\{ \{0,3,6\},\{1,5,7\},\{2,4,8\}\} \times {\bf \{\{0,3\},\{1,5\}},$ ${\bf \{2,4\}\} }$.
$$\begin{array}{lll}
\{0, 3, 6\ | \ {\bf 0, 3}\}, &\{0, 3, 6\  | \ {\bf1, 5} \},& \{0, 3, 6\ | \  {\bf2, 4}\},\\
\{1, 5, 7\ | \ {\bf 0, 3}\}, &\{1,5,7\  | \ {\bf1, 5} \},& \{1,5,7\ | \  {\bf2, 4}\},\\
\{2, 4, 8\ | \ {\bf 0, 3}\}, & \{2,4,8\  | \ {\bf1, 5} \},&\{2,4,8\ | \  {\bf2, 4}\}
\end{array}$$
In the example, sub-block $\{0,3,6\}$ appears  3 times, $\{{\bf 0,3}\}$ also appears 3 times.
The next is a non-trivial  split-block design of type $(1,1)$:
$$\begin{array}{lll}
\{0, 3, 6\ | \ {\bf 0, 3}\}, & \{0, 3, 8\  | \ {\bf1, 4} \},& \{0, 3, 7\ | \  {\bf2, 5}\},\\
\{1, 5, 7\ | \ {\bf 0, 4}\}, & \{1,5,6\  | \ {\bf1, 5} \},& \{1,5,8\ | \  {\bf2, 3}\},\\
\{2, 4, 8\ | \ {\bf 0, 5}\}, & \{2,4,7\  | \ {\bf1, 3} \},& \{2,4,6\ | \  {\bf2, 4}\}
\end{array}$$
\end{Example}
%
%%%%%%%%%%%%%%%%%%%%%%%% sec.2.2
\subsection{Some variations of split-block designs}
We modify the definition of the split-block design. In Definition \ref{split-block}, the sizes of $i$-th sub-blocks are all $k_i$.
We drop the restriction  because sub-block sizes do not need to be constant for our application to deep learning.
\begin{Definition}[Variable sub-block sizes]
Let $V_1$ and $V_2$ be mutually disjoint point sets and the block set  $\mathcal{B}$  be a collection of subsets,
each of which contains  subsets of  $V_1$ and $V_2$,  but neither subset is empty.
$$\mathcal{B} = \{ \{C_1 |  C_2\}\,  \mid \, C_1\subset V_1, C_2 \subset V_2 , \  C_1,C_2 \ne \emptyset \}.$$
For any $d_1$ and $ d_2$ points from $V_1$ and $V_2$, respectively, if there exist exactly $\lambda$ blocks containing the points, then the design is called  {\it a type $(d_1,d_2)$  split-block design with variable sub-block sizes}.
\end{Definition}

\begin{Example}
\label{VSBD}
Let $V_1=\{0,1,2,3\} , \ V_2=\{ {\bf 0,1,2}\}$. The following design is a split-block design with variable sub-block sizes of type $(2,2)$.
$$\begin{array}{lll}
\{0, 1,2\ | \ {\bf 0, 1} \}, &\{0, 1,2\  | \ {\bf 0,2} \}, & \{0, 1,2\ | \  {\bf 1,2}\}, \\
\{0,3\hspace{3.7mm}  \ | \ {\bf 0, 1}\}, & \{0,3\hspace{3.7mm} \  | \ {\bf 0,2} \},& \{0,3\hspace{3.7mm} \ | \  {\bf 1,2}\},\\
\{2,3\hspace{3.7mm} \ | \ {\bf 0, 1}\}, & \{2,3\hspace{3.7mm} \  | \ {\bf 0,2} \},& \{2,3\hspace{3.7mm} \ | \  {\bf1,2}\},\\
\{1,3\hspace{3.7mm} \ | \ {\bf 0, 1}\}, & \{1,3\hspace{3.7mm} \  | \ {\bf 0,2} \},& \{1,3\hspace{3.7mm} \ | \  {\bf1,2}\}
\end{array}$$
\end{Example}
This example satisfies the condition of split-block design with variable sub-block sizes. However,  Property \ref{regular} is not satisfied.
For example, the edge $(3,{\bf0})$ appears 6 times but $(0, {\bf0})$ appears only 4 times, that is,  this example is type $(2,2)$ but not type $(1,1)$.
Therefore, we define the split-block design with variable sub-block sizes  satisfying the  Property \ref{regular}.
%
 %%%%%%  Regular  split-block design
\begin{Definition}[Regular split-block design]
Let $(V_1, V_2\, ; \, \mathcal{B} )$  be a  type  $(d_1,d_2)$ split-block design with variable sub-block sizes.
%$$\mathcal{B} = \{ \, \{C_1 |  C_2\} \, \mid \, C_1\subset V_1, C_2 \subset V_2 \}.$$
For any $g_1$ and $g_2$ points from $V_1$ and $V_2$, respectively, $0 \le g_1 \le d_1$ and $0 \le g_2 \le  d_2$, $g_1+g_2\ge 1$,
 if there exist exactly $\lambda_{g_1,g_2}$ blocks in  $\mathcal{B}$  containing the $g_1+g_2$  points, then the design is called  a {\it regular split-block design of type $(d_1,d_2)$}.
\end{Definition}
%
%%%%%%%%
Let $V$ be a set of $v$ points, and $\mathcal{B}$ be a collection of subsets of $V$. If every $t$-subset of $V$ appears exactly $\lambda$ times in  $\mathcal{B}$, then $(V,\mathcal{B})$ is called a {\it t-wise balanced design}.
%Since block size of $t$-wise balanced design dose not have to be  the same,  usually it is not $(t-1)$-wise balanced design.
A $t$-wise balanced design with blocks in a set $K$ is not necessarily a $(t-1)$-wise balanced design.
If for any $1\le u \le t$,  $(V,\mathcal{B})$ is a $u$-wise balanced design, then it is called a {\it regular t-wise balanced design} (RtBD).
Let $(V_1,\mathcal{B}_1)$ and $(V_2,\mathcal{B}_2)$ be regular $t_1$-  and $t_2$-wise balanced designs, respectively,
then, Construction \ref{trivialD} can be generalized to   a regular split-block design
$( V_1,V_2\ ; \ \mathcal{B}_1\times \mathcal{B}_2) $.
\begin{Example}
If we add the following blocks to Example \ref{VSBD},  the combined one becomes a regular split-block design of type $(2,2)$:
$$\begin{array}{lll}
\{0\ \ \ \ | \ {\bf 0, 1} \}, &\{0\ \ \ \  | \ {\bf 0,2} \}, & \{0\ \ \ \ | \  {\bf 1,2}\}, \\
\{ 1\ \ \ \  | \ {\bf 0, 1} \}, &\{1\  \ \ \  | \ {\bf 0,2} \}, & \{ 1\  \ \ \  | \  {\bf 1,2}\}, \\
\{2\  \ \ \  | \ {\bf 0, 1} \}, &\{2\   \ \ \  | \ {\bf 0,2} \}, & \{2\  \ \ \  | \  {\bf 1,2}\}.
\end{array}$$
\end{Example}

\subsection{Related works}

We describe prior works about equivalent structure to split-block designs.
K. Ushio \cite{Ushio1981} showed a method for edge decomposition of a complete bipartite graph $K_{m,n}$ into  subgraphs isomorphic  to complete bipartite graphs $K_{a,b}$ in 1981.
Let $w(n ; k_1,k_2)$ be the number of nonnegative integer solutions $x,y$ of $n=k_1 x + k_2 y$,
where $n, k_1,k_2$ are positive integers.  We assume $n_1 \le n_2$ and $k_1 \le k_2$.
\begin{Theorem}[K. Ushio, 1981]
\  (1) When $w(n_1;k_1,k_2)=1$, a complete bipartite graph  $K_{n_1,n_2}$ has a  $K_{k_1,k_2}$ decomposition if and only if
the conditions (i)--(iv) hold.\ \
(2) When $w(n_1;k_1,k_2)\ge 2$,   a complete bipartite graph $K_{n_1,n_2}$ has a  $K_{k_1,k_2}$ decomposition if and only if the conditions (i)--(iii) hold.
\begin{description}
\item[(i)  ] $k_1k_2\, | \, n_1n_2$,
\item[(ii) ]  $ n_1 \ge k_1$ and $n_2 \ge k_2$,
\item[(iii)]  $w(n_1;k_1,k_2)\ge 1$ and $w(n_2;k_1,k_2)\ge 1$,
\item[(iv)] there exists a nonnegative integer vector $(f_1,f_2,...,f_{\beta})$ such that
$$\sum_{q=1}^{\beta} f_q = n_1 \ \mbox{and } \ k_1x_0n_2=\sum_{q=1}^{\beta} k_2y_qf_q ,$$
where $(x_0,y_0)$ is the only one solution vector of $n_1 = k_1x + k_2 y$, and $(x_q, y_q)$ for $q=1,...,\beta$ are solution vectors of $n_2=k_1x+k_2y$,

\end{description}
\end{Theorem}

D. Hoffman and M. Liatti \cite{Hoffman1995} obtained the same result  in 1995.
The decomposition of complete bipartite graph $K_{m,n}$ into  $K_{a,b}$ is equivalent to  a split-block design of type $(1,1)$ with $|V_1|=m, |V_2|=n$  and variable sub-block sizes  $\{ a , b \}$ (both of  the 1st and the 2nd sub-block sizes)  and $\lambda=1$.

In 1998, W. Martin \cite{Martin1998} defined a design similar to a split-block design called a {\it mixed $t$-design}.
\begin{Definition}[Mixed $t$-design]
Let $V_1$ and $V_2$ be the point sets of sizes $v_1, v_2$, respectively, and the block set $\mathcal{M}$ be a collection of subsets of the points consisting of $k_1$ points from $V_1$ and $k_2$ points from $V_2$:
$$ \mathcal{M} = \{ \{ C_1 \, | C_2\} \, | \, C_i\subset V_i, |C_i|=k_i , i=1,2\}.$$
For any integers $d_1,d_2$ such that $d_1+d_2 = t$,
if there exist exactly $\lambda_{d_1,d_2}$ blocks containing $d_1,d_2$ points from $V_1,V_2$, respectively,
%for any $d_1, d_2$ points such that $d_1+d_2=t$ from $V_1, V_2$,
%respectively, there exist exactly $\lambda_{d_1,d_2}$ blocks containing the points,
then the collection $\mathcal{M}$ is called a \textit{mixed $t$-design}.
% with parameters $t$-$(v_1,k_1,v_2,k_2)$.
%, where $\Lambda = (\lambda_{d_1,d_2} : d_1 + d_2 \leq t)$.
\end{Definition}
  % \begin{Definition}{(mixed $t$-design)}
  % Let $V_1$ and $V_2$ be two point sets with $|V_i| = v_i \ (i=1,2)$.
  % For fixed $k_1,k_2$ satisfying $1 \leq k_i < v_i$, let $\Omega_i$ denote the collection of all $k_i$-element subsets of $V_i$.
  % Let us refer to elements $(b_1,b_2)$ in $\Omega_1 \times \Omega_2$ as {\it dominoes}.
  % A collection $D \subseteq \Omega_1 \times \Omega_2$ of dominoes is called a {\it mixed $t$-design}
  % with parameters $t$-$(v_1,k_1,v_1,k_2,\Lambda)$
  % (where
  % $$
  % \Lambda = (\lambda_{(j_1,j_2)} : j_1 + j_2 \leq t)
  % $$
  % is an indexed family of parameters)
  % if,
  % every $(j_1,j_2)$ with $0 \leq j_i \leq k_i$ and $j_1 + j_2 \leq t$,
  % and for every pair $(S_1,S_2)$ satisfying $|S_i| = j_i$, $S_i \subseteq V_i$,
  % the number of dominoes $(b_1,b_2)$ in $D$ satisfying $S_i \subseteq b_i (i=1,2)$ is exactly $\lambda_{(j_1,j_2)}$.
  % \end{Definition}

%A mixed $t$-design is equivalent to a split-block design of type $(d_1,d_2)$ with $d_1 + d_2 = t$.
%%%%
%From the definition,  for any non-negative integers $u_1, u_2$ satisfying  $u_1+u_2 \le t$, we can say that there exist exactly $\lambda_{u_1,u_2}$ blocks containing any $u_1$ and  $u_2$ points in $V_1$ and $V_2$, respectively.
A mixed $t$-design is a split-block design of type $(d_1,d_2)$, $d_1+d_2 = t$.
%In the terminology of split-block design, the mixed $t$-design is a split block design of types $(0,t), (1,t-1), \ldots, \mbox{ and } (t,0)$.

\begin{Theorem}[W. Martin, 1998]
If there is a symmetric $2$-$(v,k,\mu)$ design,
then
there exists a mixed $2$-design with parameters  $v_1=k$, $k_1=\mu$,
$v_2=v-k$, $k_2=k-\mu$, $\lambda_{2,0} +1 = \lambda_{1,1} = \lambda_{0,2} = \mu$.
\end{Theorem}

\begin{Theorem}[W. Martin, 1998]
If there is a $3$-$(4n,2n,n-1)$ design,
then
there exists a mixed $3$-design with parameters $v_1=v_2=2n$, $k_1=k_2=n$, $\lambda_{1,1} = 2n-1$, $\lambda_{2,1} = \lambda_{1,2} = n-1$.
\end{Theorem}

If the block set of a  $t$-$(v,k,\lambda)$ design is partitionable  into classes such that
every point appears $\alpha$ times in each class, then the design is called $\alpha$-{\it resolvable block design}.  The classes are called $\alpha$-{\it resolution classes}.

K. Ozawa et al.\,\cite{ozawa2000} showed  constructions of split-block designs using $\alpha$-resolvable block designs.  The construction is basically a direct product method but
fewer copy blocks are needed.

\begin{Theorem}[K. Ozawa et al., 2000]\label{Th:bbbd}
Let $n \geq 3$ be an integer and $q$ be a prime power.
There exists a split-block design of types $(1,2)$ and $(2,1)$ with parameters
\begin{eqnarray*}
v_1 = v_2 = q^{n-1}, \quad b = \frac{q^2 (q^{n-1} -1)}{q-1}, \quad
\lambda_{11} = \frac{q^{n-1}-1}{q-1}, \quad
\lambda_{12} = \lambda_{21} = \frac{q^{n-2}-1}{q-1}.
\end{eqnarray*}
\end{Theorem}

\begin{Theorem}[K. Ozawa et al., 2000]\label{Th:bbbd2}
If there are two $\alpha_i$-resolvable $2$-$(v_i,k_i,\lambda_i)$ designs, $i=1,2$,
then there exists a split-block design of types $(1,2)$ and $(2,1)$ with parameters
$$
 \lambda_{12} = \alpha_1 d_2 \lambda_2, \quad
\lambda_{21} = \alpha_2 d_1 \lambda_1,
$$
where
$$
d_1 = \frac{\alpha_1 \mbox{lcm}(\lambda_{10}/\alpha_1,\lambda_{01}/\alpha_2)}{\lambda_{10}},\
d_2 = \frac{\alpha_2 \mbox{lcm}(\lambda_{10}/\alpha_1,\lambda_{01}/\alpha_2)}{\lambda_{01}}
$$
\end{Theorem}
M. Mishima et al. \cite{Mishima2001}  discussed {\it balanced bipartite block designs} in 2001.
This design is equivalent to a split-block design with variable sub-block sizes  but  with a constant super-block size and type $(1,1)$.

\begin{Definition}[Balanced bipartite block design]
Let $V_1$ be a set of $v_1$ points, $V_2$ be another set of $v_2$ points and $\mathcal{B}$ be a collection of $k$-subsets, called blocks (super-blocks), of $V_1 \cup V_2$.
$(V_1,V_2,\mathcal{B})$ is called a {\it balanced bipartite block design} with parameters $v_1,v_2,b,r_1,r_2,k$, $\lambda_{20},\lambda_{02},\lambda_{11}$, if
\begin{itemize}
\item[(1)] the number of replications for each point in $V_i$ is $r_i$ and any two distinct points of $V_i$ occur together in $\lambda_{20}, \lambda_{02}$ blocks, for $i=1,2$, respectively
\item[(2)] any two distinct points from different sets occur together in $\lambda_{11} $ blocks.
\end{itemize}
\end{Definition}
M. Mishima et al. showed constructions of balanced bipartite designs in the paper \cite{Mishima2001} in 2001, but they are similar to Theorems \ref{Th:bbbd} and \ref{Th:bbbd2}.

\section{Dropout designs}
\subsection{Extension of split-block designs and dropout designs}
In the previous section, we have seen designs of two layers which balance the  edges of weight.
The actual deep learning models have more than two layers.  First, we  extend it to a split-block design  having
more than two layers.

\begin{Definition}[Extended regular split-block design]
Let $V_1, V_2, \ldots, V_t$ be the mutually disjoint point sets and
$$\mathcal{B} = \{ \, \{C_1 |  C_2 | \cdots | C_t\} \,  \mid \,  C_i\subset V_i , \, C_i\ne\emptyset, \,   i=1,2,\dots,t  \}$$
be the  block set.
For any $g_1, g_2,\ldots,g_t$,  $0 \le g_i \le d_i$,  points from $V_1,V_2,\ldots, V_t$, respectively,
if there exist exactly $\lambda_{g_1,g_2,\ldots, g_t}$ blocks containing these $g_1+g_2+\cdots+g_t$  $( \ge 1 )$ points, then the design $(V_1,V_2,\ldots,V_t \, ; \, \mathcal{B})$   is called an {\it extended regular  split-block design of type $(d_1,d_2,\ldots,d_t)$}.
$t$ is said to be the {\it strength}.
\end{Definition}

Let $\{ {i_1},{i_2},\cdots ,{i_m}\}$ be a subset of  $\{ 1,2,\ldots, t\}$.
 $\mathcal{B}|_{V_{i_1}V_{i_2}\cdots V_{i_m}}$  is the set of restricted blocks ${ B}'$ from ${B} \in \mathcal{B}$
 such that  ${ B}' =\{ V_{i_1}\cup V_{i_2}\cup \cdots \cup V_{i_m}\} \cap  { B}$ for each ${ B}\in \mathcal{B}$.
 We call the block set {\it the restricted block set to} $V_{i_1}V_{i_2}\cdots V_{i_m}$.
\begin{Lemma}
Let $(V_1, V_2, \ldots, V_t\, ;\, \mathcal{B})$ be an extended regular  split-block design of type $(d_1,d_2,$ $\ldots,d_t)$.
Let $\{i_1, i_2,\ldots , i_m\}$ be a subset of $\{ 1, 2,\ldots, t\}$.
Then, $( V_{i_1},V_{i_2},\ldots , V_{i_m}  \, $;$ \, \mathcal{B}|_{V_{i_1}V_{i_2}\cdots V_{i_m}} )$ is  an extended regular  split-block design of type $(d_{i_1}, d_{i_2},\ldots,$  $d_{i_m})$.
\end{Lemma}
\Proof
Let $X_{i_j}$ be a subset of $V_{i_j}$, for $j=1,2,\ldots, m$. Since $(V_1, V_2, \ldots, V_t\, ;\, \mathcal{B})$ is an extended regular split-block design of type $(d_1,d_2,$ $\ldots,d_t)$.
The number of blocks containing
$X_{i_j}, j=1,2,..., m$ is  $\lambda^{(i_j)}$.  So  $( V_{i_1},V_{i_2},\ldots ,$ $V_{i_m}  \, ; \, \mathcal{B}|_{V_{i_1}V_{i_2}\ldots V_{i_m}} )$ is an extended regular split-block design of  type  $(d_{i_1}, d_{i_2},\ldots,$  $d_{i_m})$.
\qed

\begin{Definition}[Dropout design]
\label{defdropout}
Let $V_1, V_2, \ldots, V_n$ be the mutually disjoint point sets and
$$\mathcal{B} = \{ \, \{C_1 |  C_2 | \cdots | C_n\} \,  \mid  \,  C_i\subset V_i , ~ C_i\ne\emptyset, ~  i=1,2,\ldots,n  \}$$
the set of super-blocks.
If $(V_{i},V_{i+1},\ldots, V_{i+t-1} \, ; \, \mathcal{B}|_{V_{i}V_{i+1}\cdots V_{i+t-1}} )$
is an extended regular split-block design of type $(d_1,d_2,\ldots,d_t)$  for $i=1,2,\ldots, n-t+1$,
then $\mathcal{D}=(\, V_1, V_2, \ldots, V_n\, ; \mathcal{B})$ is called a {\it dropout design of type $(d_1,d_2,\ldots,d_t)$}.
$\lambda^{(i)}, \, i=1,2,\ldots, n-t+1$,  are the concurrence numbers of $\mathcal{D}$.
\end{Definition}

%%%%%%%%%%%%%%%%%%%% Figures 1,2
\begin{figure}[H]
\centerline {
\includegraphics[width=3.5in,clip]{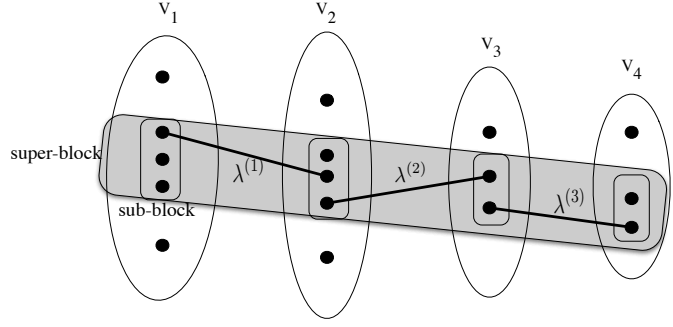}
} \caption[]{A dropout design of type (1,1)  }
\end{figure}
\begin{figure}[H]
\centerline {
\includegraphics[width=3.5in,clip]{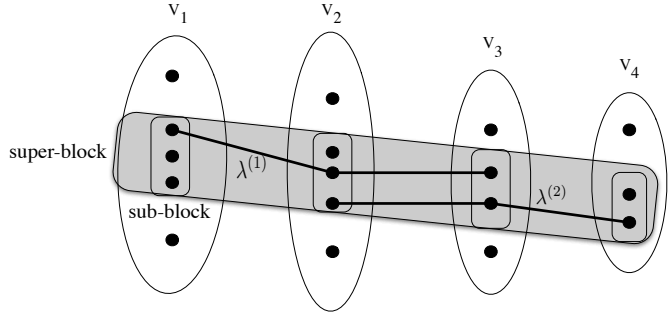}
} \caption[]{A dropout design of type (1,1,1)  }
\end{figure}

Let $B_s= \{C_{s,i} \,|\, C_{s,i+1} \,|\, \cdots \,|\, C_{s,i+t-1} \}, \, s=1,2,\ldots, b$, be the super-blocks of $\mathcal{B}|_{V_{i}V_{i+1} \cdots V_{i+t-1}}.$ Then we have the following equation:
\begin{Property}
\label{ConcNum}
$$\lambda^{(i)} \prod_{j=i}^{i+t-1} { v_j \choose d_j} =   \sum_{s=1}^b \prod_{j=i}^{i+t-1}{ |C_{s,j}| \choose  d_j  }$$
If the sizes of the $i$-th sub-blocks are all $k_i$, then the above equation can be expressed more simply as follows:
$$\lambda^{(i)} \prod_{j=i}^{i+t-1} { v_j \choose d_j} =   b \prod_{j=i}^{i+t-1}{ k_j \choose  d_j  }$$
\end{Property}
\Proof
We prove the first equation. The number of ways to choose $d_j$ points from each sub-block $C_{s,j}$ of the  block $B_s= \{C_{s,i} \,|\, C_{s,i+1} \,|\, \cdots \,|\, C_{s,i+t-1} \}$   is $\prod_{j=i}^{i+t-1}{ |C_{s,j}| \choose  d_j  }$,   and the total number for all blocks is $\sum_{s=1}^b \prod_{j=i}^{i+t-1}{ |C_{s,j}| \choose  d_j  }$.
 It is equal to the concurrence  number $\lambda^{(i)}$  times  the number of ways to choose $d_j$ points from each $V_j$ for $j=i, i+1,\ldots,i+t-1$, which is $\prod_{j=i}^{i+t-1} { v_j \choose d_j}$.
 The second equation is in the case that the block size of the  $j$-th sub-block is  $k_j$ for each $j=i,i+1,..., i+t-1$,
that is, $ |C_{s,j}| = k_j$, for any $1\le s \le b$.
\qed

The concurrence number  $\lambda_{d_1,d_2,...,d_t}$ of  \ $(V_{i},V_{i+1},..., V_{i+t-1} \, ; \, \mathcal{B}|_{V_{i}V_{i+1}\cdots V_{i+t-1}} \, )$
may  vary  for each $i=1,2,..., n-t-1$, therefore we denote it $\lambda^{(i)}$.

\begin{Example}
\label{dropout}
Let $V_1=\{0,1,2,3,4,5\}, V_2=\{ {\bf 0,1,2,3,4,5}\}, V_3=\{ {\it 0,1,2,3,4,5}\}$. The following is the dropout design of type $(1,1)$.
$$\begin{array}{lll}
\{0, 3\ | \ {\bf 0, 3} \ |\ {\it 0,3}\}, & \{0, 4\  | \ {\bf 1,5} \ |\ {\it0,5}\},& \{0, 5\ | \  {\bf2,4}\ |\ {\it0,4 }\},\\
\{1, 3\ | \ {\bf 1,4} \ |\ {\it 2,3}\}, & \{1, 4\  | \ {\bf 2,3} \ |\ {\it2,5}\},&  \{1, 5\ | \  {\bf0,5}\ |\ {\it2,4 }\},\\
\{2, 3\ | \ {\bf 2,5} \ |\ {\it 1,3}\}, & \{2, 4\  | \ {\bf 0,4} \ |\ {\it1,5}\},& \{2, 5\ | \  {\bf1,3}\ |\ {\it 1,4 }\}
\end{array}$$
\end{Example}
In the above example, the  point sets $V_1, V_2$ and  $V_3$ have $v=6$ points each,  and the sizes of sub-blocks are all $k=2$. This kind of dropout designs is easier to construct and has convenient properties. This  will be discussed  next.

%%%%%%%%     Uniform Dropout design
\subsection{Uniform dropout design}
\begin{Definition}[Uniform dropout design]
For a dropout design $\mathcal{D}=(V_1,V_2,\dots,V_n \, ;\, \mathcal{B} ) $,  if  the size of each  $V_i$ is $v$ and the size of all sub-blocks are  $k$, then  $\mathcal{D}$ is called a {\it uniform dropout design}.
\end{Definition}
\begin{Lemma}
 All concurrence numbers of  a  uniform dropout design are  the same.
\end{Lemma}
\Proof
Let  $\mathcal{D}=  (V_1,V_2,\dots,V_n \, ;\, \mathcal{B} ) $  be a uniform dropout design of type $(d_1,d_2,...,d_t)$, in which  $|V_i|=v$ for all $i=1,2,..,n$, the size of each sub-block is $k$ and the number of super-blocks is $b$.
Suppose  the concurrence numbers are $\lambda^{(i)},\  1\le i \le n-t+1$.
Consider the concurrence number $\lambda^{(1)}$ of the restricted block set $\mathcal{B}|_{V_1V_2\cdots V_t}$.
Since $v_i=v$ and $k_i=k$ for  each $i=1,2,...,n$ in  Property \ref{ConcNum},  we have
$$ \lambda^{(1)} = b \prod_{i=1}^t {k \choose d_i } / \prod_{i=1}^{t} {v \choose d_i }.$$
Similarly, we have the same equation for $\lambda^{(2)}$ of $\mathcal{B}|_{V_2V_3\cdots V_{t+1}}$.
So, $\lambda^{(1)} = \lambda^{(2)} = \cdots = \lambda^{(n-t+1)}=\lambda$.
%$$ \lambda  \prod_{i=1}^{t} {v \choose d_i }= b \prod_{i=1}^t {k \choose d_i }$$
\qed \\  \\
%%%%%%
%%%%%%%
We denote a uniform dropout design as a $(v, k, \lambda \,; n)$-{\it uniform dropout design} or  $(v, k, \lambda \,; n)$-UDD of type $(d_1,d_2,...,d_t)$.
Uniform dropout designs are very convenient for applications. First, we show a method to increase the number of layers.
\begin{Theorem}
\label{UDDtheorem}
Let $\mathcal{D}$ be a   $(v,k,\lambda\, ; t)$-UDD of type $(d_1,d_2,\dots,d_t)$.
If the design $\mathcal{D}$ is also type $(d_{2}, d_{3}, \ldots , d_t, d_{1})$, $(d_3,\ldots, d_t,d_1,d_{2})$, $\ldots$
and $(d_t,d_1,\ldots,d_{t-1})$,
then there exist a $(v,k,\lambda \,;\, n)$-UDD with the same types for any $n \ge t$.
\end{Theorem}
\Proof
Suppose that $\mathcal{D} = (V_1,\ldots,V_t; \mathcal{B})$ and $\mathcal{B} = \{ \, \{C_1 | \cdots | C_t\} \,  \mid  \,  C_i\subset V_i , ~ C_i\ne\emptyset, ~  i=1,2,\dots,t  \}$.
We extend the point sets and  block to
$$
( V_1, ..., V_t, V_1,V_2,...V_{t}\ ; \
\mathcal{B}^{'} = \{\{C_1 | \cdots | C_t | C_1 | C_2 | \ldots |C_t \} \ \}  ).
$$
Consider a consecutive $t$ point sets $V_i, V_{i+1}, ..., V_t, V_1,...V_{i-1}$.
Those restricted system forms  a  type  $(d_i, d_{i+1},...,d_t, d_1,...,d_{i-1})$ dropout design.
Since the dropout design  $\mathcal{D}$ is also type $(d_i, d_{i+1},...,d_t,$ $d_1,...,d_{i-1})$.
Likewise, it can be extended many times.
\qed

\begin{Example}
\label{UDDexample}
The following is a $(6,2,1; 2)$-UDD of type $(1,1)$.  $V_1= \{ 0,1,2,3,4,5\}, V_2=\{ {\bf 0,1,2,3,4,5}\}$.
$$
\begin{array}{ccc}
\{ 0,3 \ |\  {\bf 0,3} \} & \{0,4\  | {\ \bf 1,5}\} & \{0,5\  |\  {\bf 2,4 } \}\\
\{ 1,3\ |\   {\bf 1,4} \} & \{1,4\ |\  {\bf 2,3}\} & \{1,5\ |\  {\bf 0,5 } \}\\
\{ 2,3\ |\  {\bf 2,5} \} & \{2,4 \ |\   {\bf 0,4}\} & \{2,5\ |\   {\bf 1,3 } \}\\
\end{array}
$$
The below  is the expansion of the above to a $(6,2,1; 3)$-UDD of type $(1,1)$. $V_3 =\{ {\it 0, 1, 2, 3, }$  ${\it 4, 5}\}$.
$$
\begin{array}{ccc}
\{ 0,3 \ |\  {\bf 0,3} \ |\   {\it 0,3}\} & \{0,4 \ |\   {\bf 1,5} \ |\   {\it 0,4} \} & \{0,5\ |\   {\bf 2,4 } \ |\  {\it 0,5} \}\\
\{ 1,3 \ |\   {\bf 1,4} \ |\  {\it 1,3} \} & \{1,4 \ |\   {\bf 2,3} \ |\  {\it 1,4} \} & \{1,5 \ |\  {\bf 0,5 } \ |\  {\it 1,5} \}\\
\{ 2,3 \ |\  {\bf 2,5} \ |\   {\it 2,3} \} & \{2,4 \ |\   {\bf 0,4} \ |\   {\it 2,4} \} & \{2,5\ |\  {\bf 1,3 } \ |\   {\it 2,5} \}\\
\end{array}
$$
\end{Example}

The layer of actual deep learning models are not usually the same size.  Uniform dropout designs are
easier to construct,  but it is harder to use.  Next, we will consider  adjusting  the uniform dropout design to a more practical model.

Let $R$ be a subset of the point sets $V_1\cup V_2 \cup \cdots \cup V_n$ of a dropout design $\mathcal{D}=(V_1,V_2,...,V_n\,; \mathcal{B})$.
Consider  sub-designs of $\mathcal{D}$ whose points are reduced by  $R$. The point sets  deleted by $R$ are  $V_i' = V_i \setminus R$ for  $i=1,2...,n$.
The block set $\mathcal{B}(R)$ is the set of the modified blocks
$$\mathcal{B}(R) = \{ B\setminus R \, |\,  B\in \mathcal{B} \}, \ \
\mbox{ where }  B\setminus R = \{C_1\cap V_1' \, |\, C_1\cap V_2'  \, | \, \cdots |  C_n\cap V_n'\}. $$
Let $\mathcal{D}(R)=(V_1',V_2',...,V_n'\, ; \mathcal{B}(R))$ .
We should note that a sub-block of $B\setminus R$ can  be the empty set.  If there is an empty sub-block in a super-block of $\mathcal{B}(R)$, then the super-block cannot be a block of dropout design. Therefore, the blocks in $\mathcal{B}(R)$ having an empty sub-block are  removed.

\begin{Theorem}
Let  $\mathcal{D}=(V_1,V_2,...,V_n\,; \mathcal{B})$ be a dropout design of type $(d_1,d_2,...,d_t)$ and $R$ a subset of $V_1\cup V_2 \cup \cdots \cup V_n$.
If $R$ does not include any sub-block of $\mathcal{B}$, then the  $\mathcal{D}(R)$  is a dropout design of type $(d_1,d_2,...,d_t)$.
\end{Theorem}
\Proof
Consider $V_1,V_2,...,V_t $ as  arbitrary $t$ consecutive point sets from $V_1,V_2,...,V_n $ without loss of generality.
Suppose $\mathcal{B}(R)$ does not have empty sub-block.
Let $X_i$ be a $d_i$-subset of $V_{i}',  \, i=1,2,...,t$.
The number of blocks of $\mathcal{D}(R)$ containing  $X_i  \subset V_{i}' ,\,  i=1,2,...,t$,  is the same for any $d_i$-subsets $X_i \subset V_{ i+j-1}' $.
Thus,
$\mathcal{D}(R)$ is a dropout design which has the same type and concurrence numbers of $\mathcal{D}$.
\qed

 $\mathcal{D}(R)=(V_1',V_2',...,V_n'\, ; \mathcal{B}(R))$  is called  {\it a dropout design with deleted R}.

 \begin{Example}
\label{deletedR}
We delete $R_1 =\{0\in V_1 , {\bf 3}\in V_2 \}$ from Example \ref{dropout}.
$$\begin{array}{lllllllll}
\{ 3 &| \ {\bf 0 } & |\ {\it 0,3}\}, & \{ 4  &| \ {\bf 1,5}  &|\ {\it0,5}\},& \{ 5 &| \  {\bf2,4} &|\ {\it0,4 }\},\\
\{1, 3&| \ {\bf 1,4} &|\ {\it 2,3}\}, & \{1, 4  &| \ {\bf 2}  &|\ {\it2,5}\},&  \{1, 5 &| \  {\bf0,5} &|\ {\it2,4 }\},\\
\{2, 3&| \ {\bf 2,5} &|\ {\it 1,3}\}, & \{2, 4 &| \ {\bf 0,4}  &|\ {\it1,5}\},& \{2, 5 &| \  {\bf1} &|\ {\it 1,4 }\}
\end{array}$$
This is a dropout design with deleted $R_1$.

Next we delete $R_2=\{ {\bf 0,3} \} \subset V_2$ from  Example \ref{dropout}.
Since $\{ 0,3 \ |\ \ \emptyset  \ \ |\   {\it  0,3 }\} $  contains  an empty set, it should be removed.
Each of $0, 3\in V_1$  appear twice, but the remaining points of $V_1$ appear 3 times.
So, the design with deleted $R_2=\{ {\bf 0,3} \} $ is not a dropout design.
$$ \begin{array}{lll}
 & \{0,4 \ |\   {\bf 1,5} \ |\   {\it 0,5} \} & \{0,5\ |\   {\bf 2,4 } \ |\  {\it 0,4} \}\\
\{ 1,3 \ |\   {\bf 1,4} \ |\  {\it 2,3} \} & \{1,4 \ |\   {\bf 2 \ \ \ } \ |\  {\it 2,5} \} & \{1,5 \ |\  {\bf 5 \ \ \ } \ |\  {\it 2,4} \}\\
\{ 2,3 \ |\  {\bf 2,5} \ |\   {\it 1,3} \} & \{2,4 \ |\   {\bf 4\ \ \ } \ |\   {\it 1,5} \} & \{2,5\ |\  {\bf 1 \ \ \ } \ |\   {\it 1,4} \}\\
 \end{array}$$
\end{Example}

%%%%%%%%%%%%%%%%%%%%%%%%%%%%%%%%%%% Complementary  Design
\subsection{Complementary  dropout designs}
When we want to have a dropout design with large sub-block sizes, for instance, more than half of each $|V_i|$,
the following property of complimentary designs   is useful.

 Now we consider the set of blocks of two layers,   $\mathcal{B}|_{V_i V_j}$,  where $|i-j| \le t-1$.
 Let $X$ be a $g$-point set  in $V_i$,
 $$\mathcal{B}_i(X\, ;  V_j ) = \{ B \cap  V_j \ |\  B   \supset X  ,  B \in \mathcal{B}|_{V_i V_j} , i\ne j\}$$
 be a  set of $V_j$  part  of $\mathcal{B}|_{V_i V_j}$,  each of which contains the point set $X \subset V_i$.

\begin{Lemma}
\label{RegtwBD}
Let X be a  $w$-point set  of $V_1$, $0\le w \le d_1$.
$(V_1,V_2\, ; \mathcal{B}) $ is a dropout design of type $(d_1,d_2)$ if and only if
 $\mathcal{B}_1(X; V_2 ) $ is a regular $d_{2}$-wise balanced design for any $X\subset V_1, \ 0 \le |X| \le d_1$.
\end{Lemma}
\Proof
Suppose that $(V_1,V_2\, ; \mathcal{B}) $ is a dropout design of type $(d_1,d_2)$.
Let X be a set of $w$-points of $V_1$,  $0\le w \le d_1$.   Consider $\mathcal{B}_1(X; V_2)$.
For any  $u$-points  $Y$ of $V_2$ and any $1\le u \le d_2$,
$\mathcal{B}_1(X; V_2)$ contains $\lambda_{w,u}$ blocks, each of which includes $Y$.
So $\mathcal{B}_1(X; V_2)$ is a regular $d_2$-wise balanced design.
Conversely, suppose $\mathcal{B}_1(X; V_2)$ is a regular $d_2$-wise balanced design for any subset $X$ of $w$ points in $V_1$.
For any $w$-point $X$ in $V_1$ and $u$-point $Y $ in $V_2$,
$0 \le w \le d_1$ and $0 \le u \le d_2$ , $w+u \ge 1$,
$X$ and $Y $ simultaneously appear in the same number of blocks.
 Therefore $(V_1,V_2\ ; \mathcal{B}) $ is a dropout design of type $(d_1,d_2)$.
\qed
%%%%%%%%%%%%%%%%%%

In a  dropout design $(V_1,V_2,...,V_n \, ; \mathcal{B})$, we sometime consider a restricted system consisting of  consecutive $t$ layers
$(V_i, V_{i+1}, ..., V_{i+t-1} \, ;  \mathcal{B}|_{ V_iV_{i+1}\cdots V_{i+t-1} } )$.
Without loss of generality, we simply consider   $(V_1,V_2,...,V_t \, ; \mathcal{B}) $.
$\mathcal{B}_{12...(t-1)}(X;V_t)$ is in the case that $X$ is a subset of $V_1\cup V_2\cup \cdots \cup V_{t-1}$.

\begin{Lemma}
\label{RegtwBD2}
Let $X$ be  a set of  $w_1+w_2+\cdots+w_{t-1}\,  ( \ge 1) $ points, where each $w_i$ points are  from  $V_i$, $0 \le w_i \le d_{i}$ and $1\le i\le t-1$.
Then,     $(V_1,V_2,...,V_t \, ; \mathcal{B}) $ is a dropout design of type $(d_1,d_2,...,d_t)$  if and only if
$\mathcal{B}_{12...(t-1)}(X;V_t)$ is  a regular $d_t$-wise balanced design for any $X$.
\end{Lemma}
\Proof
Let $X_i \subset V_i$,  $|X_i| \le d_i$  and $X= X_1\cup X_2 \cup \cdots \cup X_{t-1}$.  Let $U \subset V_t,\  |U| \le d_t$.
Since $(V_1,V_2,...,V_t \ ; \mathcal{B})$ is type $(d_1, d_2,..., d_t)$ dropout design,
the number of blocks containing $X \cup U$ is a constant. Every block of the  set always includes $U$ for any $|U|$-subset of $V_t$.  So $\mathcal{B}_{12...(t-1)}(X;V_t)$ is a regular $d_t$-wise balanced design.
Conversely, for any  $U \subset V_t$ , $|U| \le d_t$, the number of block containing $U$ is a  constant. Therefor the number of blocks containing $X\cup U$ in  $(V_1,V_2,...,V_t \, ; \mathcal{B}) $
is a constant.
\qed
\begin{Lemma}
If $\mathcal{D}= (V_1,V_2, ..., V_t \, ; \mathcal{B})$ is a dropout design of type $(d_1,d_2,...,d_t)$, then  the set of the $i$-th sub-blocks  $( V_i, \mathcal{B}|_{V_i})$ is a regular $d_i$-wise balanced design for each $i =1,2,...,t$.
\end{Lemma}
\Proof
Since $\mathcal{D}$ is a dropout design of type $(d_1,d_2,...,d_t)$,  it is also of type $(0,\ldots,0,1,0,\ldots$, $0)$, which implies that every point appears exactly the same times  in  the block set $\mathcal{B}|_{V_i}$.  Similarly, $\mathcal{D}$ is also of  type
$(0,\ldots,0,u,0,\ldots,0)$, for each $u,  1\le u \le d_i$, which implies that every $u$-subset of $V_i$ appears the same times in  $\mathcal{B}|_{V_i}$. This means that  $\mathcal{B}|_{V_i}$ is a  regular $d_i$-wise balanced design.
\qed

From this lemma, we can say $( V_i , \mathcal{B}|_{V_i})$  of  a dropout design  $( V_1,V_2,...,V_n\, ; \mathcal{B}) $ is a  regular $t$-wise balanced design (RtBD)  for any $1\le i \le n$.  If the dropout design is type $(d_1,d_2,...,d_m)$, the strength $t$ of the RtBD is
 $$
 t = \left\{  \begin{array}{ll}
 \max\{d_1,d_2,\ldots,d_i\} &  \mbox{if}\ \ 1 \le i \le m-1,\\
 \max\{d_1,d_2,\ldots,d_m\} &  \mbox{if}\ \ m\le i \le n-m+1, \\
    \max\{d_{m-(n-i)},d_{m-(n-i)+1},\ldots,d_m\} & \mbox{if}\ \  n-m+2 \le i \le n.
   \\
 \end{array} \right.
 $$

Now consider a complement of a dropout design $\mathcal{D}=( V_1,V_2,...,V_n\, ; \mathcal{B})$.
We take the complement  for each sub-block of $B=\{C_1|C_2|\cdots |C_n\} \in \mathcal{B}$:
$$\bar{\mathcal{B}} =\{ \bar{B} \, |\, \mbox{ for all } B\in \mathcal{B} \} , \ \  \mbox{ where }
\bar{B} = \{ V_1\setminus C_1 \, |\,  V_2 \setminus C_2 \, |\, \cdots \, |\,  V_n\setminus C_n \}. $$
In the block set, a sub-block $V_i\setminus C_i$ may happen to be the empty set. Therefore a dropout design whose blocks never include any of $V_i$ is called  a {\it proper dropout design}.

\begin{Theorem}[Complementary dropout design]
\label{Complementary}
If $\mathcal{D}=(\, V_1, V_2, \dots, V_n\, ; \mathcal{B})$ is a proper dropout design,
then $\bar{\mathcal{D}}=(\, V_1, V_2, \ldots, V_n\, ; \bar{\mathcal{B}})$ is also a dropout design.
 $\bar{\mathcal{D}}$ is called the {\it complementary dropout design of} $\mathcal{D}$.
\end{Theorem}
 \Proof
At first, we assume  the the following results of regular $t$-wise balanced designs (RtBD).
 We also assume $(V,\mathcal{B})$ includes no duplicate blocks:
\begin{itemize}
\item[. ] Let $(V,\mathcal{B})$ be a pair of point set $V$ and collection $\mathcal{B}$ of subsets of $V$.  Then $(V,\mathcal{B})$  is a  RtBD  if and only if the complement design $(V, \bar{\mathcal{B}})$ is a RtBD,  see C. Godsil  (2010) \cite{Godsil2010}.
\item[. ]  Let $(V,\mathcal{B}_1)$ and  $(V,\mathcal{B}_2)$ be block disjoint RtBDs,
then $(V,\mathcal{B}_1 \cup \mathcal{B}_2)$ is also a RtBD.
\item[. ] Let $(V,\mathcal{B}_1)$ and  $(V,\mathcal{B}_2)$ be RtBDs such that $\mathcal{B}_2 \subset \mathcal{B}_1$, then $(V,\mathcal{B}_1 \setminus \mathcal{B}_2)$ is a RtBD.
\item[. ]  Let $(V,\mathcal{B}_1)$ and  $(V,\mathcal{B}_2)$ be RtBDs, then $(V,\mathcal{B}_1 + \mathcal{B}_2)$ is also a RtBD, where "+" is the multi-set union.
\end{itemize}
Next, we prove the theorem  in the case of two layers $V_1, V_2$.
Suppose $(V_1,V_2, \mathcal{B} )$ is a dropout design of type $(d_1,d_2)$.
Let $x_1, x_2 $ be  distinct  points of $V_1$.
 From Lemma \ref{RegtwBD},
$\mathcal{B}_1(\{x_1\};V_2)$  (also $\mathcal{B}_1(\{x_2\}; V_2)$ ) is a  regular $d_2$-wise balanced design for any $x_1$ (or $x_2$) of $V_1$.
$$
\Big(\mathcal{B}_1(\{x_1\};V_2) + \mathcal{B}_1(\{x_2\};V_2)  \Big)  \setminus  \mathcal{B}_1(\{x_1,x_2\};V_2) =\mathcal{B}_1(\{x_1\};V_2) \cup \mathcal{B}_1(\{x_2\};V_2),
$$
where $+$ operation is the multi-set union.
Since
$\mathcal{B}_1(\{x_1\};V_2)$,
$\mathcal{B}_1(\{x_2\};V_2)$ and
$\mathcal{B}_1(\{x_1,x_2\};V_2)$ are RtBDs,
$\mathcal{B}_1(\{x_1\};V_2) \cup \mathcal{B}_1(\{x_2\};V_2)$
is a RtBD and has  $x_1$ or $x_2$ in  $V_1$ part.
$$
Z=\mathcal{B}|_{V_2} \setminus \Big( \mathcal{B}_1(\{x_1\};V_2) \cup \mathcal{B}_1(\{x_2\};V_2)\Big)
$$
is also a RtBD and has neither  $x_1$ nor $x_2$ in $V_1$ part.
In the complimentary design $(V_1,V_2 ; \bar{\mathcal{B}})$,
$$ \overline{Z}  = \overline{\mathcal{B}_1}( \{x_1,x_2\}; V_2   )$$
is a RtBD. Therefore, from  Lemma \ref{RegtwBD},  the proof is completed in the case  two layers and strength 2.
Let $X \subset  V_1\cup V_2\cup \cdots \cup V_t$   such that $0 \le | X\cap V_i | \le d_i $ in a dropout design
$( V_1,V_2,...,V_t \, ;\mathcal{B})$.
From $X= \{ x_1, x_2\}$, if we continue this by induction for $X=\{x_1,x_2,x_3, \ldots\}$,
 Lemma \ref{RegtwBD2} can be applied to the complementary dropout designs.
\qed\\

\section{Constructions of dropout designs}

\subsection{Projective and affine geometries}
We begin by recalling some fundamental definitions and properties from projective and affine geometries.
Let $q$ be a prime power, $d$ a positive integer and let $V_{d+1}$ denote the $(d+1)$-dimensional vector space over a finite field of order $q$, GF$(q)$.
$(t+1)$-dimensional subspaces of $V_{d+1}$ are  called $t$-{\it flats}.
$0$-flats, $1$-flats, and $(d-1)$-flats are called {\it points},
{\it lines}, and {\it hyperplanes}, respectively.
The incidence structure of the point set and the set of $t$-flats is defined by the set theoretical inclusion relation between subspaces.
A system consisiting of all the points, all the $t$-flats of $V_d$ and their incidence structure is called {\it projective geometry}, denoted by PG$(d,q)$.
Let ${\bf x}$ be a vector of $V_{d+1}$, then $\{\alpha {\bf x} :  \alpha \in  GF(q)\setminus \{0\} \}$ is $1$-dimensional vector space ($0$-flat).  $V_{d+1}\setminus \{0\}$ is partitioned into $1$-dimensional vector spaces each of which  correspond to a point  of  PG($d,q$).
 A point of PG($d,q$) is sometime represented by a vector which is a vector in the corresponding  $1$-dimensional vector space.

The number of $t$-flats of PG$(d,q)$ is
$\begin{bmatrix}
\, d+1\, \\ t+1
\end{bmatrix}
_q$,
 where
 $
\begin{bmatrix}
\, d\, \\ t
\end{bmatrix}
_q$
is the {\it Gaussian coefficient} defined by
\begin{equation*}
 \begin{bmatrix}
\, d\, \\ t
\end{bmatrix}
_q
=
\begin{cases}
\frac{(q^d -1)(q^{d-1} -1) \cdots (q^{d-t+1}-1)}{(q^{t} -1)(q^{t-1} -1) \cdots (q-1)} &\textrm{ if } 1 \leq t \leq d, \\
1 &\textrm{ if } t = 0.
\end{cases}
 \end{equation*}
In addition, the number of $t_1$-flats through a $t_2$-flat in PG$(d,q)$ is expressed as the number of $(d-t_1 - 1)$-flats in a $(d-t_2 -1)$-flat, that is,
$\begin{bmatrix}
\, d-t_2-1\, \\ d-t_1-1
\end{bmatrix}
_q$.

By removing a hyperplane ($(d-1)$-flat)  from PG$(d,q)$,  the rest of points and all flats will produce an \textit{affine geometry} of dimension $d$, denoted by AG$(d,q)$.
The affine geometry also can be defined as follows:
The point set is the vectors of $V_d$.
Let $T$ be a $t$-dimensional vector space.
And denote the coset of $T$ as $\mathcal{C}(T)= \{ T + v \, | v\in V_d\}$.
$U \in \mathcal{C}(T)$ is called a $t$-flat.  A $1$-flat and a 2-flat are called a {\it point} and a {\it line} of affine geometry AG$(d,q)$, respectively.
A $(d-1)$-flat is called a {\it hyperplane}.
The number of $t$-flats in AG$(d,q)$ is also given by
$q^{d-t}
\begin{bmatrix}
\, d\, \\ t
\end{bmatrix}
_q$.

For  $t$-flats $T$ and $U$ of AG$(d,q)$,  $T$ and $U$ are said to be {\it parallel}  if they are in the same coset.  The set of $t$-flats in a coset is called a {\it parallel class}.
The set of $t$-flats  of AG$(d,q)$ is partitionable into parallel classes and the partition is called a {\it resolution} in the design theory.
For more details about projective and affine geometries, see T. Beth, D. Jungnickel, and H. Lenz(1999)\cite{Beth1999}, J. Hirschfeld(1998)\cite{HIRSCHFELD}.

\subsection{Constructions from orthogonal array}

%We show in this section some constructions of dropout designs  using orthogonal arrays with multi-structure.
In this section, we construct dropout designs using orthogonal arrays with multi-structure.
First, we define orthogonal arrays.

Let $S$ be a set of  $q$ symbols.
%% Let $S$ be a set of cardinality $s$.
   An {\it orthogonal array}  of strength $t$, constraints $m$ and index $\rho$
   is a $(\rho q^t \times m)$-matrix  $C$ with entries from $S$ such that
   every ordered $t$-tuple of $S$ appears exactly $\rho$ times in any $t$ columns of the matrix $C$  as a row.
%   such that each $(\rho q^t \times t)$-submatrix contains each ordered $t$-tuples
%   with entries from $S$ exactly $\rho$ times as a row.
%
 Such a matrix will be denoted by OA$_\rho (t,m,q)$.
   In the case $\rho = 1$, we write OA$(t,m,q)$.

   \begin{Lemma}\label{lem:OA}
   Let $G$ be an $n \times m$ matrix over GF$(q)$.
   If any $t$ columns of $G$ are linearly independent,
   then   the matrix whose rows are from
   $$
   \{ {\bf x}G \ : \ {\bf x} \in {\rm GF}(q)^n\}
   $$
   is  an OA$_{\rho} (t,m,q)$, $\rho = q^{n-t}$.
   \end{Lemma}

In 1997, V.  Levenshtein \cite{Levenshtein1997} defined a {\it split orthogonal array} of type $(d_1,d_2)$ with index $\rho$ which is a matrix $C$ whose row is an element of $S^{N_1+N_2}$, we simply call the matrix a ``\,{\it $C$ in $S^{N_1+N_2}$}\,''.
The matrix is partitioned into $N_1$ columns and $N_2$ columns
satisfying the condition that,
in any $d_1$ columns  in the fist $N_1$ columns  and  any $d_2$ columns in the next $N_2$ columns,
every $d_1$- and $d_2$-tuples of $S$  appears exactly $\rho$ times in the matrix $C$ as a row.
It is clear the  matrix $C$ have $\rho s^{d_1 + d_2}$ rows.

We extend a split orthogonal array by partitioning into several sub-matrices.
%%%%%%
%%%%%%
\begin{Definition}[multi-split orthogonal array]
  %Let $C$ be a matrix on GF$(q)$.
   A {\it multi-split orthogonal array of type} $(d_1,d_2,$ $..., d_t)$ and index $\rho$ is
   a matrix $C$ in $S^{N_1 + N_2 +\cdots + N_t}$
   such that every vector of $S^{d_1 + d_2+ \cdots + d_{t}}$
   occurs exactly $\rho$ times in any $d_1 + d_2 +\cdots + d_{t}$ columns,
   where each element of $S^{d_i}$  are in any  $d_i$ columns of the $N_i$ columns.
\end{Definition}

A multi-split orthogonal array can be obtained from a partition of  OA with strength $t \geq 2$.
The following result extends Lemma \ref{lem:OA}.
%obtain multi split orthogonal array as follows:

%%%%%%%%%%%%%%%%%%%%%%%%Lemma 5.3
\begin{Lemma}\label{lem:MSOA}
Let $G$ be an $s \times m$ matrix over GF$(q)$ and
partitioned into $t$ sub-matrices as follows:
$$
G = [ G_1 | G_2 | \cdots | G_t ].
$$
Suppose that every $G_i$ is an $s \times k_i$ matrix, where $\sum_{i=1}^{t}k_i = m$.
If $G$ satisfies the following conditions:\\
(1) In each $G_i$, any $d_i$ vertical vectors are linearly independent,\\
(2) any $d_1 + d_2+\cdots + d_t$ vertical vectors of $G$,  where each $d_i$ vectors are from $G_i$,  are linearly independent,\\
then the array consisting of $q^s$ vectors of
$$\{ {\bf x}G \ : \ {\bf x} \in {\rm GF}(q)^s\}$$
is a multi-split orthogonal array of type $(d_1,d_2,\ldots,d_t)$ and index $q^{s-(d_1+d_2+\cdots + d_t)}$.
\end{Lemma}
\Proof
Let $A$  be a matrix  of  $\{ {\bf x}G \ : \ {\bf x} \in {\rm GF}(q)^s\}$ and it is partitioned into $n$ sub-matrices corresponding to the columns of $G$.
From the condition of $G$, obviously, any vector from GF$(q)^{d_1 + d_2 +\cdots + d_t}$ occurs $q^{s-(d_1+d_2+\cdots + d_t)}$ times in  $A$.
\qed

Note that the matrix $G$ in Lemma \ref{lem:MSOA} is called the {\it generator matrix} with parameters $((d_1,d_2,\cdots$, $d_t),(k_1,k_2,\ldots,k_t))$.
We write $x^t$ for a $t$-tuple of all $x$, $(x,x,\ldots,x)$.

Let $A$ be a matrix  over a set of integers $S=\{0,1,\ldots,q-1\}$. Suppose  $A$  is partitioned into  $A_1, A_2,\ldots, A_n$  having  $k_1,k_2,\ldots,k_n$ columns of $A$, respectively.
Entries in every column in orthogonal array are from the same set,  therefore we have to change labels of orthogonal array  for application  to  a dropout design.
Let ${\bf c}_{ij}$ be the $j$-th column of sub-matrix $A_i$, for $1 \leq j \leq k_i$, $1 \leq i \leq n$.
$A_{i}^{'}$ is a matrix having $k_i$ columns  over the non-negative integers, $1 \leq i \leq n$, computed as follows:
$$
A^{'}_{i} = [ {\bf c}_{i1} + 0 ,  {\bf c}_{i2} + q,\ldots, {\bf c}_{ik_i} + (k_i-1)q ].
$$
Note that ${\bf c}_{ij} + x$ means the addition of $x$ to each element of ${\bf c}_{ij}$.
Here, we make the set of super-blocks $\mathcal{B}$ for a dropout design
which consists of $A_{i}^{'}$s:
$$
\mathcal{B} = \{B_i  :  i = 1,\ldots, b\}, \quad B_i = \{{\bf a}_{i1}^{'} | {\bf a}_{i2}^{'} | \cdots | {\bf a}_{in}^{'}\},
$$
where ${\bf a}_{ij}^{'}$ is the set of elements of the $i$-th row of $A_{j}^{'}$, and $b$ is the number of rows of $A$.
Let $\mathcal{D}_A = (V_1,V_2,\ldots,V_n; \mathcal{B})$ be a design with respect to a multi-split orthogonal array $A$, where $V_i$ is the set of symbols appearing in $A'_i$.
Clearly, it holds the following theorem:

\begin{Theorem}\label{Th:1^t}
Let $A$ be a matrix in GF$(q)^{k_1 + k_2,\cdots + k_n}$ .
If every consecutive $t$ sub-matricies of $A$ is a multi-split orthogonal array of type $1^t$ and index $\rho$,
then $\mathcal{D}_A =$ ($V_1,V_2,\ldots,V_n$; $\mathcal{B}$) is a dropout design of type $1^t$
with the concurrence number $\lambda = \rho$.
\end{Theorem}

\begin{Example}
Let $G$ be a $2 \times 3$ generator matrix with parameters $(1^2,(2,1))$  over  GF$(3)$ as follows:
$$
G =
\left(
\begin{array}{cc|c}
 1 & 0 & 1 \\
 1 & 2 & 2 \\
\end{array}
\right).
$$
From Lemma \ref{lem:MSOA} and Theorem \ref{Th:1^t},
$\mathcal{D}_A$ forms a dropout design of type $(1,1)$ with $\lambda = 1$, where
% as the follows:
$$
A =
\left(
\begin{array}{cc|c}
 0 & 0 & 0 \\
 1 & 2 & 2 \\
 2 & 1 & 1 \\
 1 & 0 & 1 \\
 2 & 2 & 0 \\
 0 & 1 & 2 \\
 2 & 0 & 2 \\
 0 & 2 & 1 \\
 1 & 1 & 0 \\
\end{array}
\right)
\quad \mbox{and}\
A^{'} =
\left(
\begin{array}{cc|c}
 0 & 3 & {\bf 0} \\
 1 & 5 & {\bf 2} \\
 2 & 4 & {\bf 1} \\
 1 & 3 & {\bf 1} \\
 2 & 5 &  {\bf 0} \\
 0 & 4 &  {\bf 2} \\
 2 & 3 &  {\bf 2} \\
 0 & 5 & {\bf 1} \\
 1 & 4 &  {\bf 0} \\
\end{array}
\right).
$$
\end{Example}

Let $A$ be a matrix in GF$(q)^{k_1 + k_2 + \cdots + k_n}$ such that any consecutive $t$ sub-matrices  is a multi-split orthogonal array of type $(d_1, d_2,\ldots,d_t)$ and index $\rho$.
When $(d_1,d_2,\ldots,d_t)$ is not $1^t$,
in order that $\mathcal{D}_A$  becomes a dropout design, we have to append a supplementary block set $\mathcal{B}^{\ast}$ to $\mathcal{B}$,
which satisfies the following conditions:
\\
\begin{itemize}
\item[1] Let $X_i \subset V_i,\ 2 \leq |X_i| \leq d_i$, $i = 1,\ldots,t$. Each $X_i$ is not included in any sub-block of $B|_{V_i}$, $i = 1,\ldots,t$. There exist $\rho$ blocks in $\mathcal{B}^{\ast}$ containing $X_1 \cup X_2   \cup\cdots \cup X_t$.

\item[2] Let $Y_i \subset V_i,\ 2 \leq |Y_i| \leq d_i$, $i = 1,\ldots,t$.
There does not exist any block in $\mathcal{B}^{\ast}$ containing $Y_i$
if $Y_i$ is included in a block of $\mathcal{B}$.
\end{itemize}
When $n= 2$ or $3$,
the next two theorems provide specific results about supplementary block set $\mathcal{B}^{\ast}$.

\begin{Theorem}\label{Th:(2,1)}
%Let $A$ be a $q^s \times m$ matrix on GF$(q)^{k_1+k_2}$ which is a multi-split orthogonal array of type $(2,1)$ and index $\rho$.
Let $A$ be a split orthogonal array in GF$(q)^{k_1+k_2}$ of type $(2,1)$ with index $\rho$.
Then
%If there is a multi split orthogonal array of type $(2,1)$ and index $\rho$,
%an $s \times m$ generator matrix $G$ with parameters $((2,1),(k_1, \ldots, k_n))$ in GF$(q)$ and $t \leq s$,
%Let $G$ be an $r \times m$ generator matrix with parameters $(1^t,(k_1, \ldots, k_t))$ on GF$(q)$.
%$t \leq r$とする.
$\mathcal{D}_A = (V_1,V_2\, ; \mathcal{B})$ is a dropout design of type $(2,1)$ with the concurrence number $\lambda_{2,1}= \rho$.
Note that $\lambda_{2,1}$ is the concurrence number which
is the number of blocks containing any $2$ points and $1$ point from $V_1$ and $V_2$, respectively.

\end{Theorem}

\Proof
%Using a similar way to Theorem \ref{Th:gm},
It is clear that $|V_1| = qk_1$, $|V_2| = qk_2$ and $\mathcal{D}_A$ is at least a dropout design of type $(1,1)$.
%it is obtained the matrix $A^{'}$ which is a dropout design of type $(1,1)$.
In order to  be the type $(2,1)$, we need to append some super-blocks $\mathcal{B}^{\ast}$ to $\mathcal{B}$,
%not appearing the same row in $A^{'}$.
which satisfy the following conditions:
%A collection of these superblocks have three conditions as follows:
\begin{itemize}
\item[1.] any triple $(x,y;z)$ such that
$\{x,y\}$ from $V_1$ is not contained in $\mathcal{B}$ and $z$ in $V_2$ appears in $\mathcal{B}^{\ast}$ exactly $\rho$ times,
\item[2.] any pair $\{x,y\}$ from $V_1$ appearing in $\mathcal{B}$ is not contained in $\mathcal{B}^{\ast}$.
%\item[3] the collection forms split-block design of type $(2,1)$
\end{itemize}
%Suppose $T_i$ is a partition $V_i$ into $q$ subsets, for $i=1,2$,
From the definition $A_{i}^{'}$, the symbols appearing in ${\bf c}_{i1} +0$ are $\{0,1,\ldots,q-1\}$ and ${\bf c}_{i2} +q$ is $\{q,q+1,\ldots,2q-1\}$, and so on.
Let
$$
T_1 = \{\{0,1,\ldots,q-1\},\{q,q+1,\ldots,2q-1\},\ldots,\{(k_1 -1)q, \ldots,qk_1 -1\}\}
$$and
$$T_2 ={\bf  \{\{0,1,\ldots,q-1\},\{q,q+1,\ldots,2q-1\},\ldots,\{(k_2 -1)q,\ldots,qk_2 -1\}\}}.$$
%%%%%
Any $1$st sub-block of $\mathcal{B}$ takes one by one from each set of $T_1$,
therefore any pair from a set of $T_1$ does not appear in any sub-block of $\mathcal{B}$.
Suppose $\mathcal{B}^{\ast} = T_1 \times T_2$.
Then, the design appending $\rho$ copies of $\mathcal{B}^{\ast}$ to $\mathcal{B}$ is
a dropout design of type $(2,1)$.
\qed

\begin{Theorem}
Let $A$ be a $q^s \times m$ matrix in GF$(q)^{k_1+k_2+k_3}$ such that  any successive two sub-matricies  is a multi-split orthogonal array of type $(2,1)$ and index $\rho$.
If $\rho = k_1 = k_3$,
then
$\mathcal{D}_A = (V_1,V_2,V_3;\mathcal{B})$ is a dropout design of type $(2,1)$ with the concurrence number $\lambda_{2,1}^{(1)} = \lambda_{2,1}^{(2)} = \rho$.
\end{Theorem}
\Proof
Put $T_i = \{\{qj+l : 0 \leq l \leq q-1\} : 0 \leq j \leq k_i-1\}$, $i=1,2,3$.
In a similar way to Theorem \ref{Th:(2,1)},
 the design appending $\mathcal{B}^{\ast} = T_1 \times T_2 \times T_3$ to $\mathcal{B}$ is
a dropout design of type $(2,1)$ with the concurrence number $\lambda_{2,1}^{(1)} = \lambda_{2,1}^{(2)} = \rho$,
because
$T_1 \times T_2 \times T_3$ %contains $k_1 k_2 k_3$ blocks.
 contains $T_1 \times T_2$
 $k_3$ times
and $T_2 \times T_3$ $k_1$ times, and $\rho = k_1 = k_3$.
\qed

Next, we describe some methods how to construct generator matrices of a multi-split orthogonal array using projective geometry.

\begin{Theorem} \label{th:k-11}
Let $q$ be a prime power, $d \geq 2$,  and $k_1,k_2$ be two integers such that
 $k_1 + k_2  \leq \begin{bmatrix}
\, d+1\, \\ 1
\end{bmatrix}_q (= v)$.
There exists a dropout design of type $(1,1)$ with
$$
V_1 = \{0,1,\ldots,qk_1 -1\}, \quad
V_2 = {\bf \{0,1,\ldots,qk_2 -1\}},  \quad
\lambda = q.
$$
\end{Theorem}
\Proof
The number of the points in PG$(d,q)$ is $v$,
and vector representations of any two points are linearly independent.
Then $(d+1) \times v$ matrix $G$
whose columns are vectors of  points of PG$(d,q)$
is a generator matrix with parameters $((1,1),(k_1,k_2))$, where $k_1 + k_2 = v$.
From Lemma \ref{lem:MSOA} and Theorem \ref{Th:1^t},
we have the dropout design.
\qed

\begin{Theorem}\label{th:k-cap}
Suppose $k_1,k_2,k_3$ are three integers such that
 $k_1 + k_2 + k_3  \leq k$.
If there are $k$ points no three of which are collinear in PG$(d,q)$, $d \ge 2$ and $q$ a prime power,
there exists a dropout design of type $(1,1,1)$ with
$$
V_1 = \{0,1,\ldots,qk_1 -1\}, \
V_2 = {\bf \{0,1,\ldots,qk_2 -1\}},
\
V_3 = {\it \{0,1,\ldots,qk_3 -1\}},
\
$$
where $k_1 + k_2 + k_3 \leq k$.
\end{Theorem}

The proof of the theorem is similar to Theorem \ref{th:k-11}. We omit it.
The $k$ points in Theorem \ref{th:k-cap} are studied as $k$-cap in PG$(d,q)$, see J. Hirschfeld (1998) \cite{HIRSCHFELD}.
The following is a part of known results as maximum number of $k$:

\begin{itemize}
\item When $d=2$ and $q$ is odd, there exists a $(q+1)$-cap in PG$(2,q)$% (as a conic)
\item When $d=2$ and $q$ is even, there exists a $(q+2)$-cap in PG$(2,q)$% (as a conic and its nucleus)
\item When $d=3$ and $q \neq 2$, there exists a $(q^2+1)$-cap in PG$(3,q)$% (as an elliptic)
\end{itemize}

\begin{Theorem}\label{Th:pg}
There exists a dropout design of type $(2,1)$, with
$$
V_1 = \{0,1,\ldots,q(q+1)\},
\quad
V_2 =
{\bf \{0,1,\ldots,q^3\}},
\quad
\lambda = 1.
$$
\end{Theorem}

\Proof
Let $G$ be an $3 \times (q^2+q+1)$ generator matrix in GF$(q)$ having $2$ sub-matrices $G_1,G_2$,
where $G_1$ is made from the $q+1$ points on a line $L$ in PG$(2,q)$,
and $G_2$ is from  the $q^2$ points not lie on  $L$.
Any two vectors from $G_1$ and any one vector from $G_2$
are linearly independent.
Suppose $T_1 = \{\{qi+j : 0 \leq j \leq q-1\} : 0 \leq i \leq q \}$,
$T_2 =  \{\{qi+j  : 0 \leq j \leq q-1\} : 0 \leq i \leq q^2-1\}$,
and $\mathcal{B}^{\ast} = T_1 \times T_2$.
In a similar way to Theorem \ref{Th:(2,1)},
we construct a  dropout design of type $(2,1)$.
\qed

%According to choice of the points set in Theorem \ref{Th:pg},
%we can construct a dropout design having two types which is available to Theorem 3.8.
Using Theorem \ref{UDDtheorem} and the appropriate choice of point sets in Theorem \ref{Th:pg}, we can construct the following dropout designs.

\begin{Corollary}
Let $q$ be a prime power.
There exists a $(q^2,q,1;2)$-uniform dropout design (UDD) of types $(2,1)$ and $(1,2)$.
\end{Corollary}
\Proof
Let $L$ and $L^{'}$ be distinct lines meeting at a point $p$ in PG$(2,q)$.
Take the  $q$ points from $L \setminus\{p\}$ for $G_1$, and the  $q$ points from  $L^{'} \setminus\{p\}$ for $G_2$.  Then
any two (one) from $G_1$ and any one (two) from $G_2$
are linearly independent.
 We have a generator matrix  $G=[G_1 | G_2]$ with parameters $((2,1),(q,q))$ and $((1,2),(q,q))$.
%%%
%%%
Since each of  $G_1, G_2$ over GF$(q)$ have $q$ vectors
and also using Lemma \ref{lem:MSOA} and Theorem \ref{Th:(2,1)},
the size of each  sub-block is $q$
and $|V_1| = |V_2| = q^2$.
The concurrence number is $q^{3-(2+1)} = 1$ by Lemma \ref{lem:MSOA} because
the number of rows of $G$ is $3$  and $d_1 +d_2 =3$.
\qed

\subsection{Geometrical construction}
We describe construction methods for dropout designs from the incidence structure of projective or affine  geometry.
Many geometrical structures provide $t$-designs, for example,
the points of PG$(d,q)$ together with the $t$-flats of PG$(d,q)$ as blocks form a 2-$(\begin{bmatrix}
\, d+1\, \\ 1
\end{bmatrix}_q,
\begin{bmatrix}
\, t+1\, \\ 1
\end{bmatrix}_q$,
$\begin{bmatrix}
\, d-1\, \\ t-1
\end{bmatrix}_q
)$ design and
the points of AG$(d,q)$ together with the $t$-flats of AG$(d,q)$ as blocks form a 2-$(q^d,q^t,
\begin{bmatrix}
\, d-1\, \\ t-1
\end{bmatrix}_q
)$ design (see T. Beth (1999) \cite{Beth1999}).
The following incidence structure is isomorphic to one between the set of points and the set of hyperplanes of AG$(d-t,q)$.

\begin{Lemma}\label{lem:ag}
Let $d \geq 3$ be an integer and $q$ be a prime power.
Let $T_1$ be a $t$-flat in AG$(d,q)$, $2 \leq t \leq d-1$.
Suppose that $\mathcal{C}_t(T_1)=\{T_1, \ldots, T_{q^{d-t}}\}$
 is a parallel class in AG$(d,q)$ and $\mathcal{A}^*$ is the set of hyperplanes in AG$(d,q)$ such that each hyperplane contains $q^{d-t-1}$ $t$-flats $T_{i_j}$'s
, $1 \leq i_j \leq q^{d-t}$.
Then  $(\mathcal{C}_t(T_1), \mathcal{A}^*)$ is a 2-$(q^{d-t},q^{d-t-1},\lambda)$ design with
$\lambda=\begin{bmatrix}
\, d-t-1\, \\ 1
\end{bmatrix}
_q$.
\end{Lemma}

\Proof
Let $H$ be a hyperplane  of $\mathcal{A}^*$ containing a $t$-flat $T_i$.
It is clear that $H$ includes $q^{d-t-1}$ $t$-flats of $\mathcal{C}_t(T_1)$.
Thus it holds that block size $k=q^{d-t-1}$.
The number of $(d-1)$-flats of $\mathcal{A}^*$ which contain the distinct two $t$-flats $T_i$ and $T_j$ equals to the number of $(d-1)$-flats of $\mathcal{A}^*$ containing the $(t+1)$-flat through $T_i$ and $T_j$, that is, $(q^{d-t-1}-1)/(q-1)$.
\qed

The next theorem is a generalization of  Theorem \ref{Th:bbbd}  to use $t$-flats instead of hyperplanes.

\begin{Theorem}
Let $d \geq 3$ be an integer and $q$ be a prime power.
Then there exists a $(v,k,\lambda; n)$-UDD of type $(2,1)$ and type $(1,2)$, where
 $$v=q^t, \quad k=q^{t-1},
  \quad \lambda=\frac{q^{d-2}-q^{d-t-1}}{q-1}, \quad n=q^{d-t}
 $$
for $2 \leq t \leq d-1$.
\end{Theorem}
\Proof
Let $\mathcal{D}=(T_1,T_2,\ldots,T_{q^{d-t}}; \mathcal{B})$, where $T_i$'s are $t$-flats of a parallel class $\mathcal{C}_t(T_1)$
 and $\mathcal{B}$ is the set of hyperplanes in AG$(d,q)$ any of which does not contain a $t$-flat of $\mathcal{C}_t$.
Then it is clear that $v=q^t$ and $k=q^{t-1}$ since the intersection of a $t$-flat $T_i$ and a hyperplane $H$ of $\mathcal{B}$ is a $(t-1)$-flat.
In addition, the number of blocks of $\mathcal{B}$ is equal to the number of $(d-1)$ flats in AG$(d,q)$ except $\mathcal{A}^*$ in Lemma~\ref{lem:ag}.
This yields $b=q(q^d-q^{d-t})/(q-1)$.
Consider three points $P_1$, $P_2$ of $T_i$ and $Q$ of $T_j$, $1 \leq i \neq j \leq q^{d-t}$.
The number of hyperplanes containing these three points is equal to
the number of hyperplanes through a $2$-flat in AG$(d,q)$ besides the hyperplanes of $\mathcal{A}^*$, that is,
$$\frac{q^{d-2}-1}{q-1} - \frac{q^{d-t-1}-1}{q-1} = \frac{q^{d-2}-q^{n-t-1}}{q-1}(=\lambda, \textrm{say}).$$
\qed

%%%%%%%% Theorem 5.14
\begin{Theorem}\label{pg:th1}
 Let $d \geq 2$ be an integer and $q$  a prime power.
Then there exists a $(v,k,\lambda; n)$-UDD of type $(1,1)$, where
 $$v=q^{d-1}, \quad k=1,
 %\quad \lambda_{11}=\frac{q^{n-1}-q^{n-t-1}}{q-1},
  \quad \lambda=1, \quad n=q+1.
 $$
\end{Theorem}
\Proof
Let $T$ be a $(d-2)$-flat in PG$(d,q)$.
The number of hyperplanes $H_i$ containing $T$ in PG$(d,q)$ is equals to the number of points on a line, that is, $q+1$.
Let $\mathcal{D}= (H_1^*, \ldots, H_{q+1}^*; \mathcal{B})$, where $H_i^*=H_i \setminus T$, $1 \leq i \leq q+1$ and
$\mathcal{B}$ is the set of lines such that each line is not contained in $H_i$.
Hence we have $v=q^{d-1}$.
We count the number of super-blocks of $\mathcal{B}$.
The number of lines in $H_i^*$ is equals to
\begin{equation*}
 \begin{bmatrix}
  \, d \, \\ 2
 \end{bmatrix}_q
-
\begin{bmatrix}
  \, d-1 \, \\ 2
 \end{bmatrix}_q = q^{d-2}
 \begin{bmatrix}
  \, d-1 \, \\ 1
 \end{bmatrix}_q.
\end{equation*}
Thus the number of lines in PG$(d,q)$ which is not contained in any $H_i^*$ or $T$ is given by
\begin{equation*}
 \begin{bmatrix}
  \, d+1 \, \\ 2
 \end{bmatrix}_q
- (q+1) \cdot  q^{d-2}
 \begin{bmatrix}
  \, d-1 \, \\ 1
 \end{bmatrix}_q
-
\begin{bmatrix}
  \, d-1 \, \\ 2
 \end{bmatrix}_q
=q^{2(d-1)}.
\end{equation*}
This implies that $b=q^{2(d-1)}$.
In addition, it is easily shown that each line of $\mathcal{B}$ intersects $H_i^*$ at a point.
Hence we have $k=1$.
For any two points $P$ and $Q$ in $H_i^*$ and $H_j^*$, respectively, $1 \leq i \neq j \leq q+1$,
it holds that $\lambda=1$ since there is a unique line in $\mathcal{B}$ passing through two points $P$ and $Q$.
\qed

\begin{Theorem}
 Let $d \geq 3$ be an integer and $q$ be a prime power.
Then there exists a $(v,k,\lambda; n)$-UDD of type $(2,1)$ and  $(1,2)$, where
 $$v=q^{d-1}, \quad k=q,
 %\quad \lambda_{11}=\frac{q^{n-1}-q^{n-t-1}}{q-1},
  \quad \lambda=1, \quad n=q+1.
 $$
\end{Theorem}
\Proof
We consider a $(d-2)$-flat $T$ and
$q+1$ hyperplanes $H_i$'s in PG$(d,q)$ defined in Theorem~\ref{pg:th1}.
Let $\mathcal{D}= (H_1^*, \ldots, H_{q+1}^*; \mathcal{B})$, where $H_i^*=H_i \setminus T$, $1 \leq i \leq q+1$, and
$\mathcal{B}$ is the set of planes any of which meets $T$ at a point and is not contained in $H_i$.
We count the number of super-blocks of $\mathcal{B}$.
Let $R$ be a point on $T$.
The number of planes in PG$(d,q)$ containing $R$ is equals to the number of $(d-3)$-flats in $(d-1)$ flat, that is,
$\begin{bmatrix}
  \, d \, \\ 2
 \end{bmatrix}_q$.
Thus the number of planes passing through $R$ in $\mathcal{B}$ is
%\begin{equation*}
 $\begin{bmatrix}
  \, d \, \\ 2
 \end{bmatrix}_q
-
 \begin{bmatrix}
  \, d-2 \, \\ 2
 \end{bmatrix}_q
.$
%\end{equation*}
This means that
\begin{equation*}
b = \begin{bmatrix}
  \, d-1 \, \\ 1
 \end{bmatrix}_q
\cdot
\left(
\begin{bmatrix}
  \, d \, \\ 2
 \end{bmatrix}_q
-
 \begin{bmatrix}
  \, d-2 \, \\ 2
 \end{bmatrix}_q
\right)
= \frac{q^{2d-4}(q^{d-1}-1)}{q-1}.
\end{equation*}
In addition, we can see that each plane in $\mathcal{B}$ intersects $H_i$ with a line.
Hence we have $k=q$.
For any three points $P_1$, $P_2$ of $H_i^*$ and $Q$ of $H_j^*$, $1 \leq i \neq j \leq q+1$,
it shows that $\lambda=1$ since there is a unique plane in $\mathcal{B}$ passing through these three  points $P_1$, $P_2$ and $Q$.
\qed

A \textit{spread} $\mathcal{S}$ of PG$(d,q)$ by $t$-flats is defined as a set of $t$-flats which partitions the points of PG$(d,q)$.
It is shown that there exists a spread $\mathcal{S}$ of $t$-flats of PG$(d,q)$ if and only if $d+1$ is divisible by $t+1$. (see J. Hirschfeld (1998) \cite{HIRSCHFELD})

\begin{Theorem}
 Let $d \geq 3$ and $t$ be  integers such that $d+1$ is divisible by $t+1$, and $q$ be a prime power.
Then there exists a (non-proper) dropout design of type $(2,2)$, where
 $$v= \begin{bmatrix}
 \, t+1 \, \\ 1
\end{bmatrix}_q
, \quad k=\begin{bmatrix}
 \, t \, \\ 1
\end{bmatrix}_q
\textrm{or }
\begin{bmatrix}
 \, t+1 \, \\ 1
\end{bmatrix}_q ,
 %\quad \lambda_{11}=\frac{q^{n-1}-q^{n-t-1}}{q-1},
  \quad
\lambda = \begin{bmatrix}
  \, d-3 \, \\ 1
 \end{bmatrix}_q,
%\lambda=q+1,
\quad n=(q^{d+1}-1)/(q^{t+1}-1).
 $$
\end{Theorem}
\Proof
Suppose that  $\mathcal{S}=\{T_1,\ldots,T_{n}\}$ with $n=(q^{d+1}-1)/(q^{t+1}-1)$ is a spread of $t$-flats of PG$(d,q)$.
Let $\mathcal{D}= (T_1, \ldots, T_{n}; \mathcal{B})$, where $\mathcal{B}$ is the set of hyperplanes of PG$(d,q)$.
Obviously, $b=\begin{bmatrix}
  \, d+1 \, \\ 1
 \end{bmatrix}_q.
$
Note that any hyperplane in PG$(d,q)$ intersects $T_i$ with  $(t-1)$-flat or $T_i$ itself.
For any two points $P_1$, $P_2$ of $T_i$ and any $Q_1$, $Q_2$ of $T_j$, $1 \leq i \neq j \leq n$,
the number of hyperplanes of $\mathcal{B}$ containing these four points is equal to the number of hyperplanes containing a 3-flat.
Hence we have
%\end{equation*}
%\begin{equation*}
$\lambda = \begin{bmatrix}
  \, d-3 \, \\ 1
 \end{bmatrix}_q.$
%\end{equation*}
\qed

\section{Sparse filter problem in convolutional neural networks}

The most commonly used deep neural networks is the network called Convolutional Neural Network (CNN or ConvNet). In the convolutional neural network, the learned features are convolved with the input data. This two-dimensional convolution layer makes this architecture suitable for processing 2D data such as images.
A Convolutional Neural Network (CNN) consists of a number of convolutional and subsampling (pooling) layers optionally followed by fully connected layers (multi-layer neural network). The input to a convolutional layer is an $m \times m \times r$ image (matrix $\big[x_{ij}\big]$), where $m$ is the height and width of the image and $r$ is the number of channels, e.g. an RGB image has $r = 3$.
The convolutional layer will have $n$ filters (or kernels) of size $v \times v \times q$, where $v$ is smaller than the dimension of the image and $q$ can either be the same as the number of channels $r$ or smaller and may vary for each kernel.
The size of the filters (matrix $\big[h_{ij}\big]$) gives rise to the locally connected structure which are each convolved with the image to produce $n$ feature maps (matrix $\big[u_{ij}\big] )$ of size $m-v+1$.
Here, we assume $r=q=1$. Each convolution is computed by the following way:
 $$u_{ij}=\sum_{a=1}^v\sum_{b=1}^v x_{i+a,j+b}h_{ab},  \  i,j= 1,2,...,m-v+1$$

\begin{figure}[H]
\centerline {
\includegraphics[width=3.5in,clip]{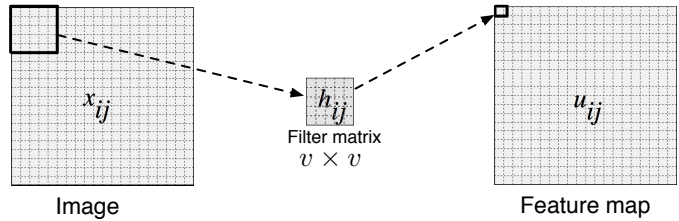}
} \caption[]{A convolution }
\end{figure}

\begin{figure}[H]
\centerline {
\includegraphics[width=3.5in,clip]{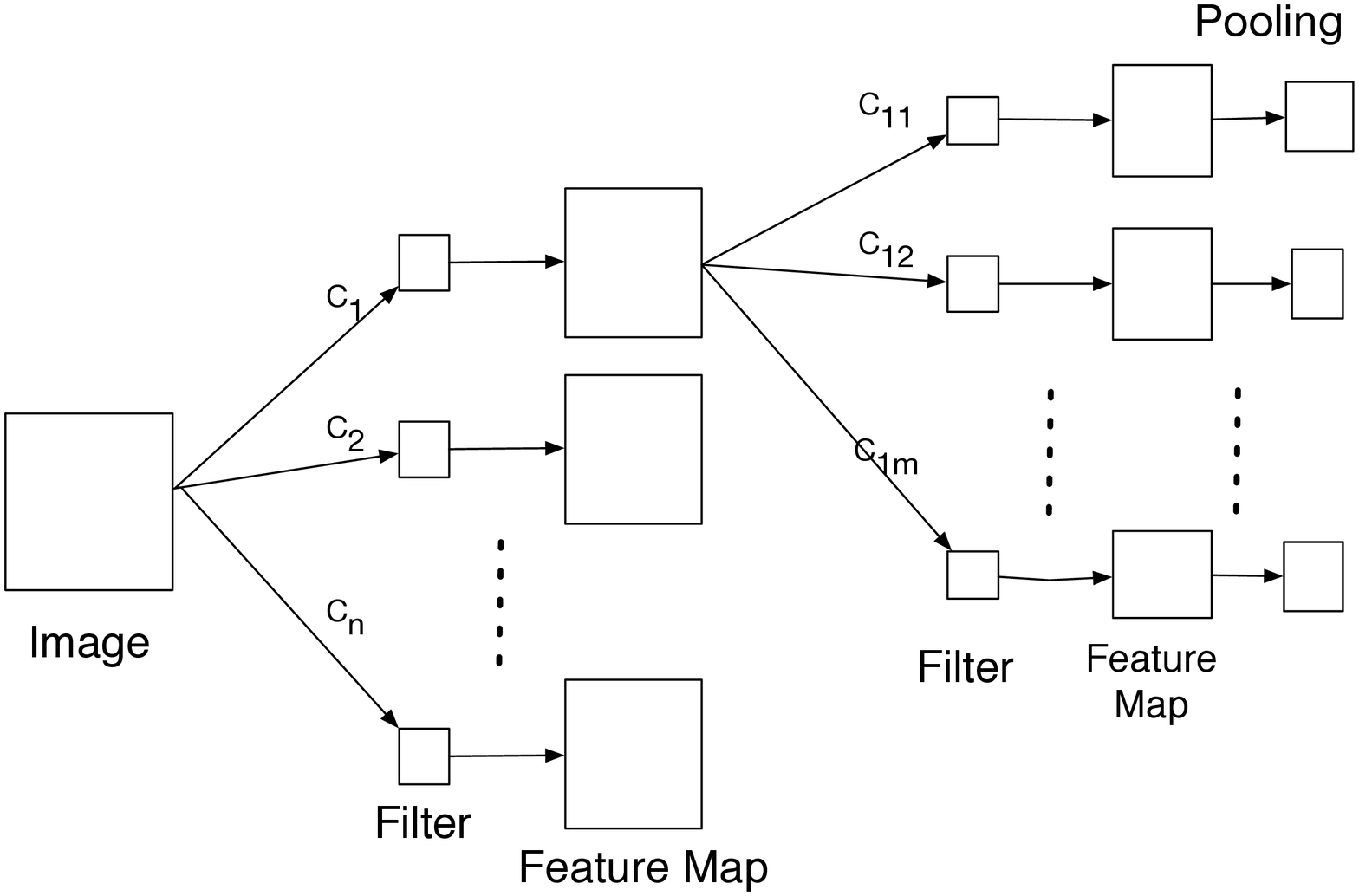}
} \caption[]{Convolutions in CNN }
\end{figure}

Initially, the filter $\big[ h_{ij}\big]$ reacts only to features that are not well understood because it contains random values, but as learning progresses, it will react strongly to features important for image recognition such as vertical lines and horizontal lines.  However, overfitting is a serious  problem in CNN too.
To prevent CNN from overfitting and to save computation time, dropout is a technique for addressing these problems.
Several methods are being discussed to solve the problem. Two kinds of dropout methods are tried, roughly classified  as follows:\\
(1) From the convolutions in the first layer $C_1, C_2, ..., C_n$, dropout several randomly at probability $p$ for each input image.
And, for the retaining units in the first layer, dropout some convolutions randomly in the second layer, e.g.  $C_{11}, C_{12},..., C_{1m}$.
N. Srivastava et al. (2014) \cite{JMLR2014} ,
W. Gao and Z. Zhou(2014) \cite{JMLR2014},
S. Changpinyo et al. (2017) \cite{DBLP2017},
Wei Wen et al. (2016) \cite{DBLP2016-2}.\\
(2) The second method is to use  sparse filters for each convolution.
A sparse filter is a $k \times k$ matrix whose randomly selected elements are preserved and the remaining elements are changed to zero.\
 N. Srivastava et al. (2014) \cite{JMLR2014},
A. Howard  et al.  (2017) \cite{DBLP2017-4},
B. Liu et al. (2015) \cite{Liu2015},
W. Wen et al. (2017) \cite{DBLP2017-2},
W. Wen et al. (2016) \cite{DBLP2016-2} and
S． Srinivas et al. (2016) \cite{DBLP2016}.

We  suggest for (1) to use a regular 2-wise balanced designs with constant block size called  2-design.
It is also called a balanced incomplete block design (BIBD), regarding existence, it has been well studied.
See  C. Colbourn and J. Dinitz (2007) \cite{Handbook}.

With respect to (2), we  propose the following $(0,1)$-matrices with  balance properties for filters.
 The integer $1$ in the matrix indicates the filter holding position and $0$ indicates the erasing position.

 \begin{Definition}
Let $H_i$ be a $v\times v$ (0,1)-matrix and $ \mathcal{B}$ be a collection of such matrices
$$ \mathcal{B} = \{ H_1, H_2, \cdots , H_b\}.$$
If  $\mathcal{B}$  satisfies the following conditions:
\begin{description}
\item[(1)]  \  Each $H_i$ has $k$ 1s in each row and column.  $ diag( H_i^T H_i) = diag(H_i H_i^T) = (k ,    k, \cdots ,  k  ) $ for each $i=1,2,..., b$,
\item[(2)]\  For each $(i,j)$ entry, the integer $1$ appears exactly $r$ times in all metrices.
$ \sum_{i=1}^b  H_i  =  r J $,
\item[(3)]\  For any distinct two rows (or columns),  the sum of their inner products for all matrices is exactly $\lambda$, that is,
  $\sum_{i=1}^b   H_i^T H_i =  \sum_{i=1}^b  H_iH_i^T  = \lambda J + (kb -\lambda)I,$
where   $J$ is  the $v\times v$  all one matrix,
\end{description}
we  call the collection of (0,1)-matrices $\mathcal{B}$  a {\it balanced filter design}.
 \end{Definition}

 In order to construct a balanced filter design, the cyclic method is useful which is  popular  in combinatorial design theory.
 Let $\mathbb{Z}_v=\{0,1,...,v-1 \}$ be  the additive group calculated with modulo $v$.  For a subset $B \subset \mathbb{Z}_v$,
  $\Delta( B) = \{ a-b \ | \ a, b \in B,  a \neq b \} $ which is a multi-set.
 Let $\mathcal{B} = \{ B_1,B_2, \dots, B_b\},  \ \ B_i \subset \mathbb{Z}_v $.
 Then the above matrix conditions are described  as follows:
 \begin{description}
 \item[(i)]\ \  \ $| B_i | = k \mbox{  for  }  i=1,2,...,b$
 \item[(ii)] \ \  $\sum_{i=1}^b  B_i = r \mathbb{Z}_v$ (which contains each element of $\mathbb{Z}_v$ $r$ times),
 \item[(iii)] \  $\sum_{i=1}^b \Delta(B_i) = \lambda ( \mathbb{Z}_v \setminus \{0\} )$,
 \end{description}
where $\sum$ is the set union as multi-set.  A collection of $k$-subsets  of $\mathbb{Z}_v$ satisfying (i) and (iii) is called {\it  difference family }.

 \begin{Example}
$$\mathcal{B}= \{ \{0,1,3\}, \{1,2,4\}, \{2,3,5\}, \{3,4,6\}, \{4,5,0\}, \{5,6,1\}, \{6,,0,2\} \}  \mbox{ in } \mathbb{Z}_7$$
satisfies the conditions (i),(ii) and (iii).
 From each  subset, e.g. $B_1=\{0,1,3\}$,  in the family, we make the following matrix.
 The first row is made that $(\{0,1,3\}+1)$-th entries  are  $1$ and other entries are $0$.
 The $j$-th row of the matrix is made by $j-1$ cyclic shifts from the first row.
 $$ H_1 = \begin{bmatrix}
 1&1&0&1&0&0&0\\
 0&1&1&0&1&0&0\\
 0&0&1&1&0&1&0\\
 0&0&0&1&1&0&1\\
 1& 0&0&0&1&1&0\\
 0&1& 0&0&0&1&1\\
 1&  0&1& 0&0&0&1
 \end{bmatrix}.$$
 For practical usage of the filter matrix, it is better to get  random permutations  $P$ and $Q $ for the rows and columns, and apply to each matrix, $ P H_1 Q ,\    P H_2 Q , ..., \ P H_b Q $. Then we can get balanced filter matrices  looks like random filters.
  $$ P  H_1 Q  = \begin{bmatrix}
0  &  1  &  1  & 0  & 0  &  1 &   0\\
1   & 0   & 1   & 0  &  1  &  0  &  0\\
0  &  1  &  0   & 1  &  1 &   0   & 0\\
 0 &   0  &  0 &   0  &  1  & 1 &  1\\
1   & 1  & 0   & 0  &  0  & 0  &  1\\
1 &  0  &  0   & 1   &0  &  1  &  0\\
0  &  0   & 1  &  1  &  0   & 0&    1
  \end{bmatrix}.$$
 \end{Example}

\section{Concluding remarks}

In this paper, we propose a new method of dropout in deep learning instead of random selection of neurons.
%The idea  is,  based on R. A. Fisher's thought,  to make balance not only the selections of nodes but also the selection of weights.
Our idea is  based on  R.A. Fisher's thought which use balance for both the selection of the nodes and the selection of the weights.
In Section 2, we investigated  combinatorial designs which  partially realize this.
In Section 3, we proposed a new combinatorial design (Definition  \ref{defdropout}) called dropout design.
Also we defined some variations of dropout designs which help to construct more practical dropout designs.
Uniform dropout designs,  dropout designs with deleted R and complementary dropout designs are defined and analyzed.
In Section 4,
estimation of weights is basically same method as regression in statistics.
Dropout in deep learning is equivalent to estimation from sparse data (including many missing data) in statistics.
It is well known in statistics what kind of sparse data is good for estimation.  We showed estimation of weights using a dropout design is based on  optimal sparse data regression.
In Section 5, we showed several constructions of dropout designs using orthogonal arrays over finite fields,
projective geometry, affine geometry and etc..
Section 6 shows sparsity problem of filters in a convolutional neural network. Instead of random sparsity, we proposed a combinatorically balanced sparse filter.

\section*{Acknowledgment}
The authors would like to thank Professor Shinji Kuriki, Osaka Prefecture University,  and Professor Ying Miao, University of Tsukuba, for valuable comments and suggestions which led to great improvement over this paper.  We would like to give  special thanks to Esther R. Lamken for many valuable comments.

%%%%%%%%%%%%%%%%%%%%%%
\bibliographystyle{abbrv}
\bibliography{dropout}

\end{document}